\documentclass[12pt,a4paper,reqno]{amsart}

\usepackage[english]{babel}
\usepackage{amssymb,latexsym,amsfonts,amsthm,upref,amsmath}
\usepackage[foot]{amsaddr}
\usepackage[margin=1in]{geometry} 
\usepackage[dvipsnames,x11names]{xcolor}
 \definecolor{myblue}{HTML}{003399}
\usepackage{hyperref}
\hypersetup{colorlinks,citecolor=myblue,filecolor=black,linkcolor=myblue,urlcolor=myblue}
\usepackage{enumerate}  
\usepackage{tikz}
\usepackage{float}
\usepackage{cleveref}
\usepackage{mathtools}
\usepackage{cite}
\usepackage{filecontents}
\makeatletter
\newcommand{\leqnomode}{\tagsleft@true}
\newcommand{\reqnomode}{\tagsleft@false}
\newcommand{\cev}[1]{\reflectbox{\ensuremath{\vec{\reflectbox{\ensuremath{#1}}}}}}
\makeatother
\newtheorem*{thm*}{Theorem}
\newtheorem*{lem*}{Lemma}
\newtheoremstyle{prim}{}{}{\normalfont}{}{\bfseries}{.}{ }{}
\newtheoremstyle{stil}{}{}{\slshape}{}{\bfseries}{.}{ }{}
\theoremstyle{stil}
\newtheorem{thm}{Theorem}[section]
\newtheoremstyle{defi}{}{}{}{}{\bfseries}{.}{ }{}
\theoremstyle{defi}
\newtheorem{defn}[thm]{Definition}
\theoremstyle{defi}
\newtheorem{rem}[thm]{Remark}
\theoremstyle{stil}
\newtheorem*{mthm*}{Main Theorem}
\newtheorem*{kor*}{Corollary}
\newtheorem{pro}[thm]{Proposition}
\theoremstyle{stil}

\theoremstyle{stil}
\newtheorem{kor}[thm]{Corollary}
\theoremstyle{prim}

\newenvironment{prf}{\noindent \textit{Proof.}}{\null\hfill$\qed$\hskip
2mm\vskip 2mm}

\newcommand{\ar}{A(\wvr{R})}
\newcommand{\arpm}{A(R^\pm)}

\newcommand{\vr}{V(\wvr{R})}

\newcommand{\aqpgl}{ \mathcal{A}_{q,p}(\widehat{\mathfrak{gl}}_N)}

\newcommand{\Yang}{ {\rm Y}(\mathfrak{gl}_{N})}

\newcommand{\Yht}{ {\rm Y} (\widehat{R})}
\newcommand{\Yhtc}{ \wtld{{\rm Y}} (\widehat{R})}

\newcommand{\Drp}{  D (\wvr{R})}

\newcommand{\tYht}{ {\rm Y} (\widehat{\mathcal{R}})}
\newcommand{\tYhtc}{ \wtld{{\rm Y}} (\widehat{\mathcal{R}})}

\newcommand{\R}{\wvr{R}}

\newcommand{\vac}{ \mathrm{\boldsymbol{1}}}

\newcommand{\gl}{\mathfrak{gl}}

\newcommand{\CC}{\mathbb{C}}

\newcommand{\ZZ}{\mathbb{Z}}

\newcommand{\Lc}{\mathcal{L}}
\newcommand{\Pc}{\mathcal{P}}

\newcommand{\Rc}{\mathcal{R}}
\newcommand{\Sc}{\mathcal{S}}

\newcommand{\Tc}{\mathcal{T}}

\newcommand{\Ec}{\mathcal{E}}

\newcommand{\Dc}{\mathcal{D}}

\newcommand{\wtld}{\widetilde}
\newcommand{\wht}{\widehat}
\newcommand{\wvr}{\overline}
\newcommand{\wndr}{\underline}

\newcommand{\ot}{\otimes}
\newcommand{\ts}{\hspace{1pt}}

\newcommand{\tr}{ {\rm tr}}
\newcommand{\sgn}{ \mathop{\rm sgn}}

\newcommand{\ndo}{\mathop{\mathrm{End}}}
\newcommand{\om}{\mathop{\mathrm{Hom}}}

\newcommand{\fand}{\quad\text{and}\quad}
\newcommand{\Fand}{\qquad\text{and}\qquad}

\newcommand{\non}{\nonumber}
\newcommand{\beq}{\begin{equation}}
\newcommand{\eeq}{\end{equation}}
\newcommand{\ben}{\begin{equation*}}
\newcommand{\een}{\end{equation*}}

\makeatletter
\def\smalloverbrace#1{\mathop{\vbox{\m@th\ialign{##\crcr\noalign{\kern3\p@}%
  \tiny\downbracefill\crcr\noalign{\kern3\p@\nointerlineskip}%
  $\hfil\displaystyle{#1}\hfil$\crcr}}}\limits}
\makeatother

\makeatletter
\def\smallunderbrace#1{\mathop{\vtop{\m@th\ialign{##\crcr
   $\hfil\displaystyle{#1}\hfil$\crcr
   \noalign{\kern3\p@\nointerlineskip}%
   \tiny\upbracefill\crcr\noalign{\kern3\p@}}}}\limits}
\makeatother

\makeatletter
\def\author@andify{%
  \nxandlist {\unskip ,\penalty-1 \space\ignorespaces}%
    {\unskip {} \@@and~}%
    {\unskip \penalty-2 \space \@@and~}%
}
\makeatother

\begin{document}

\title{Yangian deformations of $\mathcal{S}$-commutative quantum vertex algebras and Bethe subalgebras}

\author{Lucia Bagnoli}
\author{Slaven Ko\v{z}i\'{c}}
\address[L. Bagnoli and S. Ko\v{z}i\'{c}]{Department of Mathematics, Faculty of Science, University of Zagreb,  Bijeni\v{c}ka cesta 30, 10\,000 Zagreb, Croatia}
\email{lucia.bagnoli@math.hr}
\email{kslaven@math.hr}


\begin{abstract}
We construct a  new class of quantum vertex algebras associated with the   normalized Yang $R$-matrix. 
They are obtained as  Yangian deformations of certain $\mathcal{S}$-commutative quantum vertex algebras
and  their $\mathcal{S}$-locality takes the form of a single   $RTT$-relation.
 We establish some preliminary results on their representation theory and then further investigate   their braiding map. 
In particular, we show that its fixed points   are closely related with  Bethe subalgebras   in the Yangian   quantization of  the Poisson algebra $\mathcal{O}(\mathfrak{gl}_N((z^{-1})))$, which were recently introduced by Krylov and Rybnikov. 
Finally, we extend this construction of commutative families   to the case of  trigonometric $R$-matrix of type $A$.
\end{abstract}

\maketitle

\allowdisplaybreaks


\section{Introduction}\label{intro}
\setcounter{equation}{0}
\numberwithin{equation}{section}
The notion of quantum vertex algebra was introduced by Etingof and Kazhdan \cite{EK} motivated by deformed chiral algebras of E. Frenkel and Reshetikhin \cite{FR2}. It features the $\Sc$-locality property, a quantum analogue of the ordinary locality axiom for vertex algebras.
 In the first examples of quantum vertex algebras, which were given in \cite{EK} and associated with  rational, trigonometric and elliptic   $R$-matrices of type $A$, the $\Sc$-locality takes the form of the quantum current commutation relation of Reshetikhin and Semenov-Tian-Shansky \cite{RS}.
Suppressing the arguments, the relation can be written as
$$
\Lc_{13}\ts R_{12}^{-1}\ts \Lc_{23}\ts R_{12}
\sim
R_{21}^{-1}\ts\Lc_{23}\ts R_{21} \ts\Lc_{13},
$$
where $R_{ij}$ denotes the $R$-matrix and $\Lc_{ij}$ the matrix of formal power series of generators of the underlying algebra $V$, both applied on the    factors $i$ and $j$ of the triple tensor product $\ndo\CC^N\ot\ndo\CC^N\ot V$; see \cite{EK,RS} for more information. 

In this paper, we consider the problem of constructing quantum vertex algebras such that their $\Sc$-locality is in the form of the  $RLL$-relation, which is also known as $FRT$-relation due to   the pioneer work \cite{FRT} of Faddeev, Reshetikhin and Takhtajan,
\beq\label{rtttfrom}
R_{12}\ts \Lc_{13}\ts \Lc_{23} \sim \Lc_{23}\ts \Lc_{13}\ts R_{12}.
\eeq
More specifically, we are interested in  structures  which are not $\Sc$-commutative, so that their vertex operator map possesses   nontrivial singular part; cf. De Sole,   Gardini and Kac \cite[Sect. 4]{DGK}.
Our research is motivated by the fact that there exist numerous classes of quantum algebras  defined by the $RLL$-relations such that the matrix $\Lc=\Lc(u)$ consisting  of formal power series of the algebra generators possesses infinitely many positive and negative powers of the variable $u$. 
For example, in the case of rational and trigonometric $R$-matrices, such algebras naturally occur in the study of antidominantly shifted Yangians; see, e.g., the papers by  
 Frassek,  Pestun and Tsymbaliuk \cite{FPT} and Krylov and Rybnikov \cite{KR}. On the other hand, in the case of elliptic $R$-matrix of the eight-vertex model, the defining relations of such a form produce  the  elliptic quantum algebra   $\mathcal{A}_{q,p}(\widehat{\mathfrak{sl}}_2)$ which goes back to 
 Foda,   Iohara,  Jimbo,   Kedem,   Miwa and Yan \cite{FIJKMY}; see also the more recent generalization to $\mathcal{A}_{q,p}(\widehat{\mathfrak{gl}}_n)$ by Frappat, Issing and Ragoucy \cite{FIR}. Hopefully,   better understanding of the aforementioned problem will lead to  the development of vertex algebraic framework   for    such algebras and their representations, as we partially  demonstrate in this paper using  the setting of \cite{KR}. 

The paper is organized as follows. In Section \ref{section01}, we  set the notation. Next, in Section \ref{sec0201}, we introduce the algebra $\Drp$ 
associated with the normalized Yang $R$-matrix $\wvr{R}(u)$. Its definition resembles that of the double Yangian for the general linear Lie algebra $\mathfrak{gl}_N$. It possesses two families of generators organized into matrices $L^\pm (u)\in\ndo\CC^N\ot\Drp[[u^{\pm 1}]]$, where
one family generates the Yangian for $\mathfrak{gl}_N$ and another one a certain algebra $\ar$  (which differs from the dual Yangian).  The later subalgebra is naturally equipped   with the structure of $\Sc$-commutative quantum vertex algebra as well as with the structure of $\Drp$-module. In Section \ref{sec0202}, we use the  Yangian for $\mathfrak{gl}_N$ to deform the former, $\Sc$-commutative structure, thus getting the quantum vertex algebra $\vr$.
Its   vertex operator map is expressed in terms of the operator series 
in $\ndo\CC^N\ot\ndo\vr [[u^{\pm 1}]]$,
$$
L(u) = \left( L^+(u)^t \ts  L^-(u)^t\right)^t,
$$
where the superscript $t$ denotes the matrix transposition in $\ndo\CC^N$, and its $\Sc$-locality property   takes the form \eqref{rtttfrom} due to $L(u)$ satisfying the $RLL$-relation
$$
  \R_{12}  ( u-v ) \ts  L_1(u)\ts L_2 (v)
=
 L_2 (v)\ts  L_1(u)\ts \R_{12}  ( -v+u).
$$

Our next goal is to apply this vertex algebraic framework to the algebra $\Yht$, recently introduced in \cite{KR},  which can be regarded as 
the Yangian quantization of the Poisson algebra $\mathcal{O}(\mathfrak{gl}_N((z^{-1})))$. In order to do so, we investigate the following two aspects of the quantum vertex algebra $\vr$.
In Section  \ref{sec08}, we derive some elementary results on its representation theory which, as we demonstrate in Section \ref{section06}, imply that $\vr$-modules are naturally equipped with the structure of $\Yht$-module. Next, in Section
\ref{sec07}, we find families of fixed points of the braiding map of $\vr$. One of the main results in \cite{KR} is the construction of Bethe subalgebras, associated with the anti-symmetrizer,  of a certain completion of $\Yht$, which  then give rise to   Bethe subalgebras in the corresponding antidominantly shifted Yangians.  In Section \ref{section06}, we show that, in the sense of the aforementioned representation-theoretic connection, the   generators of these subalgebras come from the fixed points of the braiding of $\vr$. Moreover, we prove that some of these generators belong to the center of the completed algebra $\Yht$
and we extend the original construction to the case of symmetrizer.
Finally, in Section \ref{zadnjisection}, motivated by the possible applications to the antidominantly shifted $RTT$ quantum affine algebras of $\mathfrak{gl}_N$ \cite{FPT}, we extend   some  results from \cite{KR} and Section \ref{section06} to the case of trigonometric $R$-matrix in type $A$. More specifically, we obtain explicit formulae for commutative families and a family of central elements in the trigonometric analogue of the completed algebra $\Yht$.

\section{Preliminaries}\label{section01}

In this section, we introduce some notational conventions which are used throughout the paper. Let $N\geqslant 2$ be an integer and $h$ a formal parameter. Suppose that $V$ is a module for the commutative ring $\CC[[h]]$ of all formal power series in $h$ with complex coefficients. Recall that $V$ is said to be   topologically free  if it is isomorphic to $V_0[[h]]$ for some complex vector space $V_0$ or, equivalently, if it is torsion free, separated and complete with respect to the $h$-adic topology; cf. \cite[Ch. XVI]{Kas}. Suppose $V$ is a topologically free $\CC[[h]]$-module. We shall write $V((z))_h$ for the $\CC[[h]]$-module of all formal series
$$
a(z)=\sum_{r\in\ZZ} a_r z^{-r-1}\in V[[z^{\pm 1}]]
\quad\text{such that}\quad \lim_{r\to\infty}a_r=0, 
$$
where the limit is taken with respect to the $h$-adic topology. Note that   $V((z))_h$ coincides with the $h$-adic completion of $V((z))$, so that we have $V((z))_h=V_0((z))[[h]]$ for some complex vector space $V_0$.

Throughout the paper, we employ the usual   expansion convention, where the expressions of the form $(a_1 u_1+\ldots +a_n u_n)^r$ with $r<0$ and nonzero $a_1,\ldots ,a_n\in\CC$ are expanded in nonnegative powers of the variables $u_2,\ldots ,u_n$. 
For example, we have
$$
(u_1-u_2)^{-1}=\sum_{r\geqslant 0}\frac{u_2^r}{u_1^{r+1}}\neq
-\sum_{r\geqslant 0}\frac{u_1^r}{u_2^{r+1}}=(-u_2+u_1)^{-1}.
$$
The same applies to the case when   the last variable $u_n$ is replaced by the  parameter $h$.

Let $R(u)$ be a   $R$-matrix with coefficients in $\ndo\CC^N\ot \ndo\CC^N$ which satisfies   the   Yang--Baxter equation 
\beq\label{ybe}
R_{12}(u)\ts R_{13}(u+v)\ts R_{23}(v)=R_{23}(v)\ts R_{13}(u+v)\ts R_{12}(u).
\eeq
Note that the above equality is given over $(\ndo\CC^N)^{\ot 3}$
and 
we use the subscripts $1,2,3$  to indicate the corresponding tensor factors so that, e.g., we have
$R_{12}(u)=R(u)\ot 1$.

Let $n$ be a positive integer. 
Define the  $R$-matrix product  $ R_{[n]}(u)=R_{[n]}(u_1 , \ldots ,u_n )$ with entries in
$(\ndo\CC^N )^{\ot n}  $, which  depends  on the   
    variables $u = (u_1 , \ldots ,u_n )$,   
 by
$$
 R_{[n]}(u)= \prod_{i=1,\dots,n-1}^{\longrightarrow} 
\prod_{j=i+1,\ldots,n}^{\longrightarrow}  R_{ij}(  u_i -u_j ).
$$
The arrows in the above definition  indicate the order of factors, so that, e.g., we have
$$
R_{[4]}(u)=R_{12}\ts R_{13}\ts R_{14}\ts R_{23}\ts R_{24}\ts R_{34},\quad\text{where} \quad R_{ij}=R_{ij}(u_i-u_j).
$$

Let $m$ be a positive integer and $a\in\CC$. Introduce the products depending on the variable $z$, families of variables $u = (u_1 , \ldots ,u_n )$ and $v = (v_1 , \ldots ,v_m )$ and the formal parameter $h$ with entries in  
$(\ndo\CC^N )^{\ot n} \ot (\ndo\CC^N )^{\ot m}$
 by
\begin{align}\label{oppositeof}
 R_{nm}^{12}(z+u-v+ah)= \prod_{i=1,\dots,n}^{\longrightarrow} 
\prod_{j=n+1,\ldots,n+m}^{\longleftarrow}   R_{ij} ( z+u_i-v_{j-n}+ah ). 
\end{align}
Also,  omitting the variable $z$  we have 
\beq\label{oppositeof2f}
R_{nm}^{12}( u-v+ah )= \prod_{i=1,\dots,n}^{\longrightarrow} 
\prod_{j=n+1,\ldots,n+m}^{\longleftarrow}   R_{ij} ( u_i-v_{j-n}+ah ).
\eeq
We shall write bar on the top of the numbers in the superscript to indicate that the $R$-matrices in the given product    come  in the order    opposite to \eqref{oppositeof}, i.e.
\begin{align}
&R_{nm}^{\bar{\scriptstyle 1} 2}( x )= \prod_{i=1,\dots,n}^{\longleftarrow} 
\prod_{j=n+1,\ldots,n+m}^{\longleftarrow}   R_{ij}
, \quad
R_{nm}^{1\bar{\scriptstyle 2}}( x )= \prod_{i=1,\dots,n}^{\longrightarrow} 
\prod_{j=n+1,\ldots,n+m}^{\longrightarrow}   R_{ij}, \label{jos2s}\\
 &R_{nm}^{\bar{\scriptstyle 1}\bar{\scriptstyle 2}}(x)= \prod_{i=1,\dots,n}^{\longleftarrow} 
\prod_{j=n+1,\ldots,n+m}^{\longrightarrow}   R_{ij},
\quad\text{where}\quad 
R_{ij} =R_{ij} ( z+u_i-v_{j-n}+ah ) \label{jos}
\end{align}
and $x$ stands for $z+u-v+ah$. We use the analogous notation in the case of  $x= u-v+ah$, with the original arguments in \eqref{jos} now replaced by $ u_i-v_{j-n}+ah$. For example, if $n=m=2$ and $R_{ij}=R_{ij}(u_i-v_{j-2}+ah)$ we have
$$
R_{22}^{12}(u-v+ah)=R_{14}R_{13}R_{24}R_{23}\fand
R_{22}^{\bar{\scriptstyle 1} 2}(u-v+ah)=R_{24}R_{23}R_{14}R_{13}.
$$

\section{On certain algebras associated with the Yang $R$-matrix}\label{sec0201}

In this section, we associate certain algebras with the   Yang $R$-matrix  
\beq\label{betheyang}
R (u)= I -P \ts h\ts u^{-1}\in\ndo\CC^N\ot \ndo\CC^N[ h/u ],
\eeq
where $I$ and $P$ are the identity and the  permutation operator on $\CC^N\ot\CC^N$.
First, let us normalize $R(u)$ as follows,
\begin{align}
 R^{\pm}(u)=u (u\mp h)^{-1} R(u)  \in \ndo\CC^N\ot \ndo\CC^N[[ h/u ]].\label{yang}
\end{align}
In order to consider both normalizations simultaneously, we denote by $\wvr{R}(u)$ any of the $R$-matrices $R^+(u)$ and $R^-(u)$.
Naturally, both $R(u)$ and $\wvr{R}(u)$   satisfy the Yang--Baxter equation \eqref{ybe}. Moreover, the $R$-matrix
$\wvr{R}(u)= R^{\pm}(u)$ possesses the unitarity property
\beq\label{uni}
\wvr{R}(u)\ts \wvr{R}(-u) =1.
\eeq

Let $\ar$ be the topologically free  associative algebra over the ring $\CC[[h]]$  generated by the elements $l_{ij}^{(-r)}$, where $r=1,2\ldots  $ and $i,j=1,\ldots ,N$,   subject to the defining relations 
\beq\label{rtt}
\R ( -u+v )\ts  L_1^+(u)\ts L_2^+ (v)
=  L_2^+ (v)\ts L_1^+(u)\ts\R ( u-v  ).
\eeq
The matrix of generators $L^+(u) $ is defined by
\beq\label{ofthfrm}
L^+(u) =\sum_{i,j=1}^N e_{ij}\ot l_{ij}^+ (u),\quad\text{where}\quad l_{ij}^+(u)= \sum_{r=1}^{\infty}l_{ij}^{(-r)}u^{r-1} 
\eeq
and $e_{ij}$ are matrix units.
In general, we use subscripts to indicate a copy of the matrix of the form \eqref{ofthfrm} in the tensor product algebra $(\ndo\CC^N)^{\ot n}\ot \ar[[u]]$ so that, e.g., 
$$
L_m^+(u) =\sum_{i,j=1}^N 1^{\ot (m-1)}\ot e_{ij}\ot 1^{\ot (N-m)}\ot l_{ij}^+ (u).
$$
In particular, we have $n=2$ and $m=1,2$ in the defining relation \eqref{rtt}.
It is worth noting that, by the unitarity property \eqref{uni},   the relation \eqref{rtt} is equivalent to
\beq\label{rtt98}
L_1^+(u)\ts L_2^+ (v)
= \R (  u-v )\ts  L_2^+ (v)\ts L_1^+(u)\ts\R ( u-v  ).
\eeq

We now construct certain representations of $\ar=A(R^\pm)$. Introduce the algebras  
\begin{align*}
&\mathcal{P}^+ =\CC[x_1,\ldots, x_N,y_1,\ldots ,y_N, t_r : r=1,2,\ldots][[h]],\\
&\mathcal{P}^- = \Lambda(x_1,\ldots, x_N)\ot \Lambda(y_1,\ldots, y_N)\ot \CC[  t_r  : r=1,2,\ldots][[h]],
\end{align*}
where $\Lambda(x_1,\ldots,x_N)$ stands for the exterior algebra in generators  $x_1,\ldots,x_N$.
 By closely examining the defining relations \eqref{rtt} one easily verifies the next proposition.

\begin{pro}\label{repje}
The assignments
$$
l_{ij}^{(-r)}\mapsto x_i \ts y_j\ts  t_r,\qquad\text{where}\quad i,j=1,\ldots , N,\,r=1,2,\ldots   ,
$$
define a structure of $A(R^{\pm})$-module on $\Pc^\pm$. Moreover, the images of the monomials
$$
l_{i_1 j_1}^{(-r_1)}\ldots l_{i_m j_m}^{(-r_m)},\quad\text{where}\quad
m\in\ZZ_{\geqslant 0},\, i_1,\ldots ,i_m,\,j_1,\ldots,j_m=1,\ldots ,N,\,1\leqslant r_1\leqslant\ldots\leqslant r_m 
$$
are such that $i_1\leqslant\ldots\leqslant i_m$ (resp. $i_1<\ldots< i_m$) and
$j_1\leqslant\ldots\leqslant j_m$ (resp. $j_1<\ldots< j_m$), form a (topological) basis of the image of $A(R^{+})$ (resp. $A(R^{-})$) under this representation.
\end{pro}

 In the next definition, we extend the algebra $\ar$ by the  Yangian $\text{Y}(\gl_N)$ of the general linear Lie algebra $\gl_N$.
Let $\Drp$ be the $h$-adically  complete topological algebra over   $\CC[[h]]$ generated by the elements $l_{ij}^{(r)}$, where $r\in\ZZ$, $r\neq 0$ and $i,j=1,\ldots ,N$, subject to the defining relations
\begin{gather}
\wvr{R} (-u+v)\ts  L^+_1(u) \ts  L^+_2(v)=
 L^+_2(v)\ts  L^+_1(u) \ts \wvr{R} (u-v),\label{DRR1}\\
R (u-v)\ts  L^-_1(u) \ts  L^-_2(v)=
 L^-_2(v)\ts  L^-_1(u) \ts R (u-v),\label{DRR2}\\
   L^-_1(u) \ts  L^+_2(v)=
	\left(\wvr{R} (-u+v)^{t_1}\ts 
 L^+_2(v)\ts  L^-_1(u) ^{t_1}\right)^{t_1},\label{DRR3}
\end{gather}
where 
  $t\colon\ndo\CC^N\to\ndo\CC^N$ denotes the matrix transposition $e_{ij}\mapsto e_{ji}$
	with subscript $1$ indicating its application on the first tensor copy of $\ndo\CC^N \ot \ndo\CC^N \ot \Drp$.
The matrices of generators   $ L^\pm(u)$ are given by \eqref{ofthfrm} and by
$$
 L^{-}(u)=\sum_{i,j=1}^N e_{ij}\ot l_{ij}^- (u),\quad\text{where}\quad
l_{ij}^-(u)=\delta_{ij}+h \sum_{r=1}^{\infty}l_{ij}^{(r)}u^{-r}
.
$$

Observe that the subalgebra of   $\Drp$ generated by all $l_{ij}^{(r)}$ with $i,j=1,\ldots ,N$ and $r=1,2,\ldots$ is exactly the Yangian $\Yang$ for   $\gl_N$. Indeed, it is clear from the defining relations and the Poincar\'{e}--Birkhoff--Witt theorem  for the Yangian  \cite{O} (see also \cite[Ch. 1]{M} for more details and references)   that there exists an algebra epimorphism   $\Drp\to \Yang$ which annihilates all $l_{ij}^{(-r)}$ for $i,j=1,\ldots ,N$ and $r=1,2,\ldots .$ As for the subalgebra generated by all $l_{ij}^{(-r)}$ with $i,j=1,\ldots ,N$ and $r=1,2,\ldots ,$ the next proposition implies that it coincides with $\ar$.

\begin{pro}\label{pro21}
There exists a unique structure of $\Drp$-module over $\ar$ such that 
for all $n\geqslant 0$  the action of the generator matrices $L^\pm (u)$ of  $\Drp$ on
$$
 L_1^+(v_1)\ldots L_n^+(v_n)  \in  (\ndo\CC^N)^{\ot n}\ot \ar [[v_1,\ldots ,v_n]]
$$
is given by
\begin{align}
&L^+_0(u )\ts  L_1^+(v_1)\ldots L_n^+(v_n)  
=
  L_0^+(u)
 L_1^+(v_1) \ldots L_n^+(v_n),\label{lemma21e}\\
&L^-_0(u )\ts  L_1^+(v_1)\ldots L_n^+(v_n)  
=
 \R_{0n}( u  -v_n )^{-1}\ldots  \R_{01}( u  -v_1 )^{-1}\ts
 L_1^+(v_1) \ldots L_n^+(v_n).\label{lemma21}
\end{align}
Moreover, the action of $L^-(u)$ belongs to 
$ \ndo\CC^N \ot \om(\ar,\ar[u^{-1}]_h )$.
\end{pro}

\begin{prf}
It is clear   that the identity \eqref{lemma21e} with $n\geqslant 0$ defines an operator series $L^+(u)$ on $\ar$.
To prove that the  series $ L^- (u )$ is well-defined  by \eqref{lemma21}, it suffices to check that it preserves the ideal of defining relations \eqref{rtt} for the algebra   $\ar$. Let $n\geqslant 2$ be an integer and $v=(v_1,\ldots ,v_n)$ a family of variables. For any $i=1,\ldots ,n-1$ consider the image of 
$ \R_{i\ts i+1 }( -v_i +v_{i+1} )L_{1}^+(v_1)\ldots L_{n}^+(v_n) $, which corresponds to the left-hand side of \eqref{rtt}, under  $ L^-_0(u )$.  By \eqref{lemma21}, it is equal  to the expression
$$
\R_{i\ts i+1 }( -v_i +v_{i+1} ) \ts \R_{0n}^{-1}\ldots \R_{01}^{-1} 
\ts L_{1}^+(v_1)\ldots \ts L_{n}^+(v_n) \quad\text{for}\quad\R_{0i}^{-1}=\R_{0i}( u  -v_i )^{-1}.
$$
Due to unitarity \eqref{uni}, we can use the Yang--Baxter equation \eqref{ybe} to rewrite it as
$$
 \R^{(i)}\ts \R_{i\ts i+1 }( -v_i +v_{i+1} )\ts   L_{1}^+(v_1)\ldots \ts L_{n}^+(v_n) ,
$$
where $\R^{(i)}=\R_{0n}^{-1}\ldots \R_{0\ts i+2}^{-1} \ts \R_{0 i}^{-1}\ts \R_{0\ts i+1}^{-1}
\ts \R_{0\ts i-1}^{-1}\ldots
\R_{01}^{-1}$.
However, by   \eqref{rtt} this equals
\begin{align*}
\R^{(i )} 
\ts L_{1}^+(v_1)\ldots \ts L_{i-1}^+(v_{i-1})\ts L_{i+1}^+(v_{i+1})\ts L_{i}^+(v_{i})\ts L_{i+2}^+(v_{i+2}) \ldots
\ts L_{n}^+(v_n) \ts \R_{i\ts i+1} (v_i -v_{i+1} ).
\end{align*}
Finally,    observe that by \eqref{lemma21} the above expression coincides with the image of
$$
 L_{1}^+(v_1)\ldots \ts L_{i-1}^+(v_{i-1})\ts L_{i+1}^+(v_{i+1})\ts L_{i}^+(v_{i})\ts L_{i+2}^+(v_{i+2}) \ldots
\ts L_{n}^+(v_n) \ts \R_{i\ts i+1} (v_i -v_{i+1} ),
$$
   which corresponds to the right-hand side of \eqref{rtt},
under $L^-_0(u)$. Thus, we proved that the operator series $L^-(u)$ is well-defined by \eqref{lemma21}. Moreover, it is clear from the form of the $R$-matrix \eqref{yang} and the unitarity property \eqref{uni} that $L^-(u)$ belongs to 
$ \ndo\CC^N \ot \om(\ar,\ar[u^{-1}]_h )$. 

To finish the proof, it remains to check that the operator  series defined by \eqref{lemma21e} and 
\eqref{lemma21} satisfy the defining relations \eqref{DRR1}--\eqref{DRR3}
of  $\Drp$. However, this is verified by a straightforward calculation which relies on the Yang--Baxter equation \eqref{ybe}.
\end{prf}

 Motivated by the last assertion of Proposition \ref{pro21}, we introduce the notion of restricted $\Drp$-module as follows. A $\Drp$-module $W$ is said to be restricted if it is a topologically free $\CC[[h]]$-module such that the action of
$L^-(u)$ belongs to 
$ \ndo\CC^N \ot \om(W,W[u^{-1}]_h )$. 
Throughout the rest of this section, we assume that $W$ is a   restricted $\Drp$-module.
Define the   operator series on $W$ by
\beq\label{preslikavanjeel}
L(u) = \left( L^+(u)^t \ts  L^-(u)^t\right)^t.
\eeq
  Note  that $L(u)$ is well-defined and, furthermore, belongs to 
$ \ndo\CC^N \ot \om(W,W((u))_h )$, due to $W$ being restricted. Moreover,  a simple calculation relying on the defining relations of $\Drp$ shows that $ L (u)$ satisfies the $RLL$-relation
\beq\label{rtt2}
  \R  ( u-v ) \ts  L_1(u)\ts L_2 (v)
=
 L_2 (v)\ts  L_1(u)\ts \R  ( -v+u). 
\eeq
Finally,  formulae \eqref{lemma21e} and \eqref{lemma21} imply that for  $W=\ar$    its action is given by
\beq\label{lemma21l}
L_0(u )\ts  L_1^+(v_1)\ldots L_n^+(v_n)  
=
 \R_{0n}( u  -v_n )^{-1}\ldots  \R_{01}( u  -v_1 )^{-1} 
 L_0^+(u)L_1^+(v_1) \ldots L_n^+(v_n).
\eeq

Our next goal is to generalize       \eqref{rtt2} and \eqref{lemma21l}. 
For any $n\geqslant 1$ introduce a family of variables $u=(u_1,\ldots ,u_n)$ and
consider the   operators on $(\ndo\CC^N)^{\ot n}\ot W$ given by
\beq\label{Tn}
L^{\pm}_{[n]}(u)=L_1^{\pm}(u_1)\ldots L_n^{\pm}(u_n)\fand  L_{[n]}(u)=\R_{[n]}(u_1,\ldots ,u_n)\ts  L_{1}(u_1)\ldots L_{n}(u_n).
\eeq
Note that $L_{[n]}(u)$ can be also expressed in terms of operator series $L^{\pm}(u)$ as
\beq\label{altexpr}
L_{[n]} ( u )=
\left( L_1^+( u_1)^{t_1}\ldots L_n^+( u_n)^{t_n}
\ts
 L_n^-( u_n)^{t_n}\ldots L_1^-( u_1)^{t_1}\right)^{t_1\ldots t_n}.
\eeq
The next proposition is proved by a direct calculation relying on \eqref{ybe}, \eqref{rtt2} and \eqref{lemma21l}. 
\begin{pro} \label{ppro12}
Let $W$ be a restricted  $\Drp$-module.
For any integers $n,m\geqslant 1$, the families of variables $u = (u_1 , \ldots ,u_n )$, $v = (v_1 , \ldots ,v_m )$ and the single variable $z$ we have
\begin{align}
&L_{[n]}(u)\in  (\ndo\CC^N)^{\ot n} \ot \om(W,W((u_1,\ldots ,u_n))_h ),\label{ghwqcd5}
\\
 &\R_{nm}^{\bar{\scriptstyle 1}\bar{\scriptstyle 2}}( u-v )   L_{[n]}^{13}(u)\ts L_{[m]}^{23} (v)
=
L_{[m]}^{23} (v)\ts L_{[n]}^{13}(u)\ts \R_{nm}^{ 12}(  -v+u  )\label{rtt2gen}.
\end{align}
Moreover, if $W=\ar$ we have
\begin{align}
& L_{[n]}^{13}(z+u)\ts L^{+23}_{[m]}(v)\vac
=
\R_{nm}^{\bar{\scriptstyle 1}\bar{\scriptstyle 2}}(z+u-v)^{-1}  
L_{[n+m]}^+ (z+u,v)\vac,\label{lemma21gen}
\end{align}
where $z+u$    denotes the $n$-tuple
$(z+u_1,\ldots ,z+u_n)$ and $\vac$ is the unit in $\ar$. In particular, we have
$ 
L_{[n]} ( u ) \vac
=
L_{[n ]}^+ ( u)\vac
$ in $ \ar$.
\end{pro}
It is worth noting that, due to \eqref{ghwqcd5} with $W=\ar$,  the expression
$$
L_{[n]}(z+u)=L_{[n]}(v_1,\ldots ,v_n)\big|_{v_1=z+u_1\ldots ,v_n=z+u_n},
$$
which appears in \eqref{lemma21gen}, is a well-defined element of
$$
(\ndo\CC^N)^{\ot n}\ot\om(\ar,\ar((z))_h[[u_1,\ldots,u_n]]).
$$

\section{Quantum vertex algebra \texorpdfstring{$\vr$}{V(R)}}\label{sec0202}

In this section, we equip $\ar$ with the structure of quantum vertex algebra.
From now on, the tensor products over $\CC[[h]]$ are understood as $h$-adically completed.
For readers' convenience, we recall the definition of quantum vertex algebra   \cite[Sect. 1.4.1]{EK}.

\begin{defn}\label{qvoa}
A  quantum vertex algebra  is a quadruple $(V,Y,\vac,\Sc)$ which satisfies the following axioms:
\begin{enumerate}[(1)]
\item  $V$ is a topologically free $\mathbb{C}[[h]]$-module.
\item $Y=Y(z)=Y(\cdot ,z
)$ is the vertex operator map, a $\mathbb{C}[[h]]$-module map
\begin{align}
Y \colon V\ot V&\to V((z))_h,\label{trunc89}\\
u\ot v&\mapsto Y(z)(u\ot v)=Y(u,z)v=\sum_{r\in\mathbb{Z}} u_r v \ts z^{-r-1},\non
\end{align}
which satisfies the   weak associativity:
for any $u,v,w\in V$ and $k \in\mathbb{Z}_{\geqslant 0}$
there exists $r\in\mathbb{Z}_{\geqslant 0}$
such that
\begin{align}
&\left((z_0 +z_2)^r\ts Y(u,z_0 +z_2)Y(v,z_2)\ts w \right.\non\\
&\left.\quad - (z_0 +z_2)^r\ts Y\big(Y(u,z_0)v,z_2\big)\ts w\right) 
\in  h^k   V[[z_0^{\pm 1},z_2^{\pm 1}]].\label{associativity}
\end{align}
\item $\vac$ is the   vacuum vector, a distinct element of $V$  which satisfies
\beq\label{v1}
Y(\vac ,z)v=v,\quad Y(v,z)\vac \in V[[z]]\fand \lim_{z\to 0} Y(v,z)\ts\vac =v\quad\text{for all }v\in V.
\eeq
\item $\Sc=\Sc(z)=1+\mathcal{O}(h)$ is the  braiding, a $\mathbb{C}[[h]]$-module map
$ V\otimes V\to V\otimes V\otimes\mathbb{C}((z))[[h]]$ which satisfies the shift condition
\begin{align}
&[\Dc\otimes 1, \mathcal{S}(z)]=-\frac{d}{dz}\mathcal{S}(z)\quad \text{for}\quad \Dc\in\ndo V\text{   defined by }\Dc v=v_{-2}\vac,\label{s1}\\
\intertext{the  Yang--Baxter equation}
&\mathcal{S}_{12}(z_1)\ts\mathcal{S}_{13}(z_1+z_2)\ts\mathcal{S}_{23}(z_2)
=\mathcal{S}_{23}(z_2)\ts\mathcal{S}_{13}(z_1+z_2)\ts\mathcal{S}_{12}(z_1),\label{s2}\\
\intertext{the  unitarity condition}
&\mathcal{S}_{21}(z)=\mathcal{S}^{-1}(-z),\label{s3}
\end{align}
the  $\mathcal{S}$-locality:
for any $u,v\in V$ and $k\in\mathbb{Z}_{\geqslant 0}$ there exists
$r\in\mathbb{Z}_{\geqslant 0}$ such that  
\begin{align}
&\left((z_1-z_2)^{r}\ts Y(z_1)\big(1\otimes Y(z_2)\big)\big(\mathcal{S}(z_1 -z_2)(u\otimes v)\otimes w\big)\right.
\nonumber\\
&\left.\quad-(z_1-z_2)^{r}\ts Y(z_2)\big(1\otimes Y(z_1)\big)(v\otimes u\otimes w)\right) \,
\in\,  h^k V[[z_1^{\pm 1},z_2^{\pm 1}]] \quad\text{for all }w\in V,\label{locality}
\end{align}
and the  hexagon identity
\begin{align}\label{hexagon}
\Sc(z_1)\left(Y(z_2)\ot 1\right) =\left(Y(z_2)\ot 1\right)\Sc_{23}(z_1)\ts \Sc_{13}(z_1+z_2).
\end{align}
\end{enumerate}
\end{defn}

In the following theorem, we associate  with the $R$-matrix $\R(u)=R^\pm (u)$ a  quantum vertex algebra. It  is given over the $\CC[[h]]$-module of $\ar$, which we denote by $\vr$    to indicate that we now   regard  its   quantum vertex algebra structure.

\begin{thm}\label{mainn}
There exists a unique   quantum vertex algebra structure on $\vr=\ar$ such that the vertex operator map   is given by
\beq\label{Ymap}
Y(L^+_{[n]}(u)\vac,z)=  L_{[n]}(z+u),
\eeq
the vacuum vector is the unit $\vac $
and the braiding     is defined by  
\beq\label{Smap}
\mathcal{S}(z)\big(  L_{[n]}^{+13}(u)  L_{[m]}^{+24}(v) 
  \vac^{\ot 2}    \big) 
 =\R_{nm}^{\bar{\scriptstyle 1}\bar{\scriptstyle 2}}( z+u-v )\ts L_{[n]}^{+13}(u) \ts  L_{[m]}^{+24}(v) \ts \R_{nm}^{ 12}( z+u-v )^{-1}
 (\vac^{\ot 2} ) .
\eeq
\end{thm}

\begin{prf}
To prove that the vertex operator map $Y(\cdot, z)$ and the braiding $\Sc (z)$ are well-defined  by \eqref{Ymap} and \eqref{Smap}, it suffices to verify that they preserve the ideal of defining relations for the algebra $\ar$. This follows by a straightforward calculation
which is carried out similarly  to the proof of  Proposition \ref{pro21}.  The calculation relies   on the Yang--Baxter equation \eqref{ybe}, defining relations \eqref{rtt} and    relation \eqref{rtt2gen}. 
Next, we observe that the image of the vertex operator map  belongs to $ \vr((z))_h $
due to \eqref{ghwqcd5}.

Our next goal is to verify weak associativity \eqref{associativity}. Let $n,m,k\geqslant 1$ be arbitrary integers and $u=(u_1,\ldots ,u_n)$, $v=(v_1,\ldots ,v_m)$, $w=(w_1,\ldots ,w_k)$ families of variables. By applying $Y(z_0 +z_2)(1\ot Y(z_2))$, which corresponds to the first summand in \eqref{associativity}, on 
\beq\label{arg}
L_{[n]}^{+14}(u)\ts L_{[m]}^{+25}(v)\ts L_{[k]}^{+36}(w)(\vac\ot\vac\ot\vac),
\eeq 
we get
\beq\label{assoc1h}
 L_{[n]}^{14}(z_0+z_2+u)\ts L_{[m]}^{24}(z_2+v)\ts L_{[k]}^{+34}(w) \vac .
\eeq
Let us turn to the second summand. Applying $Y(z_2)(Y(z_0)\ot 1)$ on  \eqref{arg} we obtain
$$
Y(\ts Y(\ts L_{[n]}^{+14}(u)\vac,z_0)\ts L_{[m]}^{+24}(v)\vac,z_2)\ts L_{[k]}^{+34}(w)\vac
=
Y( L_{[n]}^{14}(z_0+u)\ts L_{[m]}^{+24}(v)\vac,z_2)\ts L_{[k]}^{+34}(w)\vac,
$$
where the equality follows from  the definition \eqref{Ymap} of the vertex operator map. By employing \eqref{lemma21gen} we rewrite the right-hand side as
\beq\label{gft}
Y(\ts\R_{nm}^{\bar{\scriptstyle 1}\bar{\scriptstyle 2}}( z_0+u-v )^{-1} \ts 
L_{[n]}^{+14} (z_0+u)L_{[m]}^{+24}(v)\vac,z_2)\ts L_{[k]}^{+34}(w)\vac.
\eeq
As the vertex operator map and the $R$-matrices are applied on   different tensor factors, their actions commute. Thus, using \eqref{Ymap} again, we find that \eqref{gft} equals
\begin{gather}
\R_{nm}^{\bar{\scriptstyle 1}\bar{\scriptstyle 2}}( z_0+u-v )^{-1}\ts 
 L_{[n+m]}^{12\ts 4} (z_2+z_0+u,z_2+v) \ts  L_{[k]}^{+34}(w)\vac,\quad\text{where}\label{arg2}\\
(z_2+z_0+u,z_2+v)=(z_2+z_0+u_1,\ldots ,z_2+z_0+u_n,z_2+v_1,\ldots ,z_2+v_m).\label{zedva}
\end{gather}
However, by employing the Yang--Baxter equation \eqref{ybe} one   easily obtains the equality
$$
 \R_{nm}^{\bar{\scriptstyle 1}\bar{\scriptstyle 2}}( z_0+u-v )^{-1}\ts 
\R_{[n+m]}^{12} (  z_0+u,v  )
=   \R_{[n]}^1 ( u )\ts \R_{[m]}^2 ( v),
$$
where $(z_0+u,v)$ is obtained from \eqref{zedva} by omitting the variable  $z_2$.
Hence, \eqref{arg2}   equals
\beq\label{arg3}
   L_{[n]}^{14}( z_2+z_0+u)\ts L_{[m]}^{24}(z_2+v)\ts L_{[k]}^{+34}(w)\vac.
\eeq
 By using \eqref{lemma21gen} one checks that the expression in \eqref{arg3}  is equal  to
\begin{align}
&\R_{mk}^{\bar{2}\bar{3}}( z_2+v-w )^{-1}\ts
\R_{nk}^{\bar{1}\bar{3}}( z+u-w )^{-1}\ts
\R_{nm}^{\bar{1}\bar{2}}( z+u-z_2-v )^{-1}
\non
\\
&\times 
L_{[n]}^{+14}(z+u)\ts 
L_{[m]}^{+24}(z_2+v)\ts
L_{[k]}^{+34}(w)\vac\label{arg4}
\end{align}
for $z=z_2+z_0$. On the other hand, by comparing \eqref{assoc1h}  and \eqref{arg3}, we see that the analogous  expression for \eqref{assoc1h} is obtained  by setting $z=z_0+z_2$ in \eqref{arg4}. 
Finally, fix nonnegative integers $l,l_1,\ldots, l_{n+m+k}$. 
Consider the coefficients of all monomials 
\beq\label{monomials}
u_1^{p_1}\ldots u_n^{p_n} v_1^{p_{n+1}}\ldots v_{m}^{p_{n+m}}w_1^{p_{n+m+1}}\ldots
w_k^{p_{n+m+k}},\text{ }  1\leqslant p_j\leqslant l_j,\text{ }j=1,\ldots ,n+m+k,
\eeq
in \eqref{arg4}. Clearly, they possess finitely many negative powers of the variable $z$ modulo $h^l$. Therefore, multiplying the coefficients by $z^r$ for a sufficiently large positive integer $r$ we obtain a power series which possesses only  nonnegative powers of the variable $z$ modulo $h^l$. This implies that the coefficients of monomials \eqref{monomials} in the product of \eqref{assoc1h} and  $(z_0+z_2)^r$ and in the product of \eqref{arg3} and $(z_0+z_2)^r$ coincide modulo $h^l$, so that weak associativity \eqref{associativity} holds.

The axioms \eqref{v1}  concerning the vacuum vector  follow from   \eqref{Ymap} with the use of Proposition \ref{ppro12}. In particular, it is worth noting that Proposition \ref{ppro12} implies the identity
\beq\label{eq317}
Y(L_{[n]}^{+}(u)\vac, z)\vac  = L_{[n]}^{+}(z+u)\vac
\eeq
for all $n$.
The proof of shift condition \eqref{s1} goes in  parallel with the   corresponding part of the proof of \cite[Thm. 2.2]{BJK}.
More specifically, it
relies on the equalities
 $$
\Dc\vac=0\Fand
\Dc\ts L_{[n]}^+(u)\vac=\textstyle\left(\sum_{l=1}^n \frac{\partial}{\partial u_l}\right)L_{[n]}^+(u)\vac\quad\text{for}\quad u=(u_1,\ldots ,u_n),\, n\geqslant 1,
$$
which are found by computing the coefficient of   $z$ in \eqref{eq317}.  

The requirements \eqref{s2} and   \eqref{s3} imposed  on the braiding are verified by a direct calculation which employs the corresponding properties of the $R$-matrix, the Yang--Baxter equation \eqref{ybe} and the unitarity \eqref{uni}.

Let us prove the $\Sc$-locality \eqref{locality}. By applying $Y(z_2)(1\ot Y(z_1))$, which corresponds to the second summand in \eqref{locality}, to $L_{[m]}^{+23}(v)L_{[n]}^{+14}(u)(\vac\ot \vac)$, we get
\beq\label{locc1}
 L_{[m]}^{23}(z_2+v)\ts
 L_{[n]}^{13}(z_1+u).
\eeq
On the other hand, by applying $Y(z_1)(1\ot Y(z_2))\Sc(z_1 -z_2)$, which corresponds to the first summand in \eqref{locality}, on $L_{[n]}^{+13}(u)L_{[m]}^{+24}(v)(\vac\ot \vac)$, we obtain
\begin{align}
\R_{nm}^{\bar{1}\bar{2}}( z_1-z_2+u-v )\ts
 L_{[n]}^{13}(z_1+u)\ts
 L_{[m]}^{23}(z_2+v)\ts
\R_{nm}^{12}( z_1-z_2+u-v )^{-1}.\label{locc2}
\end{align}
Fix nonnegative integers $l,l_1,\ldots, l_{n+m+k}$. Consider the coefficients of all monomials \eqref{monomials} in   $ \R_{nm}^{ 12}(  z+u-v )^{-1}$. 
They possess finitely many negative powers of the variable $z$ modulo $h^l$. Suppose that $z^{-r}$ for some $r>0$ is the lowest power of $z$  modulo $h^l$ in the expression. Then the coefficients of all monomials \eqref{monomials} in the product of $(z_1-z_2)^r$ and \eqref{locc2}  coincide with the corresponding coefficients in
\begin{align}
&\big( \R_{nm}^{\bar{1}\bar{2}}( z_1-z_2+u-v )\ts
 L_{[n]}^{13}(z_1+u)\ts
 L_{[m]}^{23}(z_2+v)\big)\non \\
&\times 
\big((z_1-z_2)^r 
\R_{nm}^{ 12}( -z_2+z_1+u-v )^{-1}\big)\mod h^l.\label{idalje}
\end{align}
Finally, by employing \eqref{rtt2gen} we rewrite this as
\beq\label{locc4}
(z_1-z_2)^r \ts 
 L_{[m]}^{23}(z_2+v) \ts
 L_{[n]}^{13}(z_1+u)\mod h^l.
\eeq
Clearly,  \eqref{locc4} is equal to the product of $(z_1-z_2)^r$ and \eqref{locc1}. 
Therefore, if we  multiply  \eqref{locc1} and \eqref{locc2} by $(z_1 -z_2)^r$, we obtain expressions whose  coefficients with respect to  the monomials \eqref{monomials}
 coincide   modulo $h^l$, thus proving the $\Sc$-locality.

It remains to prove the hexagon identity. Throughout its proof, we use the notation
$$
\R_{nm}^{\bar{1}\bar{2}}=\R_{nm}^{\bar{1}\bar{2}}( z_2+u-v ),\quad
\R_{nk}^{\bar{1}\bar{3}}=\R_{nk}^{\bar{1}\bar{3}}( z_1+z_2+u-w ),\quad
\R_{mk}^{\bar{2}\bar{3}}=\R_{mk}^{\bar{2}\bar{3}}( z_1+v-w ).
$$
 Consider  the left-hand side of  \eqref{hexagon}. First, by applying $Y(z_2)\ot 1$     on \eqref{arg}
we get
$$
 L_{[n]}^{14}(z_2 +u)\ts L_{[m]}^{+24}(v)\ts L_{[k]}^{+35}(w)(\vac\ot\vac).$$
By \eqref{lemma21gen} this is equal to
\beq\label{hexp1}
\big(\R_{nm}^{\bar{1}\bar{2}}\big)^{-1} 
L_{[n]}^{+14}(z_2+u)\ts 
L_{[m]}^{+24}(v)\ts 
L_{[k]}^{+35}(w)(\vac\ot\vac).
\eeq
Finally, by applying the braiding $\Sc(z_1)$   to \eqref{hexp1} we get
\begin{align}
&\big(\R_{nm}^{\bar{1}\bar{2}}\big)^{-1} \ts
\R_{mk}^{\bar{2}\bar{3}} \ts
\R_{nk}^{\bar{1}\bar{3}} \ts
L_{[n]}^{+14}(z_2+u)\ts L_{[m]}^{+24}(v)\ts 
L_{[k]}^{+35}(w)\ts
\big(\R_{mk}^{ 23}\big)^{-1} 
\big(\R_{nk}^{ 13}\big)^{-1}
(\vac\ot\vac).\label{hexp2}
\end{align}

Consider the right-hand side of \eqref{hexagon}. Applying $\Sc_{13}(z_1 + z_2)$ on \eqref{arg} we obtain
$$
\R_{nk}^{\bar{1}\bar{3}} \ts
L_{[n]}^{+14}(u)\ts 
L_{[m]}^{+25}(v)\ts 
L_{[k]}^{+36}(w)\ts
\big(\R_{nk}^{ 13}\big)^{-1}(\vac\ot\vac\ot \vac).
$$
Next, we apply $\Sc_{23}(z_1)$, thus getting
\begin{align*}
&\R_{nk}^{\bar{1}\bar{3}} \ts
\R_{mk}^{\bar{2}\bar{3}} \ts
L_{[n]}^{+14}(u)\ts 
L_{[m]}^{+25}(v)\ts L_{[k]}^{+36}(w)\ts
\big(\R_{mk}^{ 23}\big)^{-1} 
\big(\R_{nk}^{ 13}\big)^{-1}(\vac\ot\vac\ot\vac).
\end{align*}
Finally, by applying $Y(z_2)\ot 1$ and then using \eqref{lemma21gen}    we find
\begin{align}
&\R_{nk}^{\bar{1}\bar{3}} \ts
\R_{mk}^{\bar{2}\bar{3}} \ts
\big(\R_{nm}^{\bar{1}\bar{2}}\big)^{-1} 
L_{[n]}^{+14}(z_2+u)\ts 
L_{[m]}^{+24}(v)\ts 
L_{[k]}^{+35}(w)\ts
\big(\R_{mk}^{ 23}\big)^{-1} 
\big(\R_{nk}^{ 13}\big)^{-1}(\vac\ot\vac).\label{hexp4}
\end{align}
It remains to observe that the identity
$$
\R_{nm}^{\bar{1}\bar{2}} \ts
\R_{nk}^{\bar{1}\bar{3}} \ts
\R_{mk}^{\bar{2}\bar{3}}  
= 
\R_{mk}^{\bar{2}\bar{3}} \ts
\R_{nk}^{\bar{1}\bar{3}} \ts
\R_{nm}^{\bar{1}\bar{2}} 
,
$$
which is a simple generalization of the Yang--Baxter equation \eqref{ybe}, implies   the equality of the expressions in \eqref{hexp2} and \eqref{hexp4}, so the hexagon identity   follows.
\end{prf}

It is evident from the proof that, due to  relation \eqref{rtt2}, the $\Sc$-locality property  \eqref{locality}    possesses the $RLL$-form as in \eqref{rtttfrom}.   Furthermore, observe that the vertex operator map has a nontrivial singular part, which comes from the action of the Yangian $\Yang$ on $\vr$. However, at the classical limit $h\to 0$ this action vanishes, thus giving rise to a commutative vertex algebra. 
Recall \eqref{altexpr}. 
In view of the alternative expression 
$$
Y(L_{[n]}^+(u),z)=
\left( L_1^+( z+u_1)^{t_1}\ldots L_n^+( z+u_n)^{t_n}
\ts
 L_n^-( z+u_n)^{t_n}\ldots L_1^-( z+u_1)^{t_1}\right)^{t_1\ldots t_n}
$$
 for the vertex operator map \eqref{Ymap}, one can, roughly speaking, regard $\vr$ as a Yangian deformation of the quantum  vertex algebra structure over $\ar$ which is given by
$$
Y(L_{[n]}^+(u),z)=
  L_1^+( z+u_1) \ldots L_n^+( z+u_n) .
$$
Clearly, such a vertex operator map is $\Sc$-commutative, i.e. it satisfies the $\Sc$-locality \eqref{locality} with $r=0$ and the braiding defined by
$$
\mathcal{S}(z)\big(  L_{[n]}^{+13}(u)  L_{[m]}^{+24}(v) 
  \vac^{\ot 2}    \big) 
 =\R_{nm}^{12}( -z-u+v )\ts L_{[n]}^{+13}(u) \ts  L_{[m]}^{+24}(v) \ts \R_{nm}^{ 12}( z+u-v )^{-1}
 (\vac^{\ot 2} ) .
 $$

\section{On \texorpdfstring{$\vr$}{V(R)}-modules}\label{sec08}

In this section, we investigate a connection  between $\vr$-modules and restricted $\Drp$-modules.
For reader's convenience, we start by recalling the definition of module for       quantum vertex algebra   \cite[Def. 2.23]{Li}.

\begin{defn}\label{qvoamodule}
 Let $(V,Y,\vac ,\Sc)$  be a  quantum vertex algebra. A     $V$-module  is a pair $(W,Y_W)$, where $W$ is a topologically free $\CC[[h]]$-module
and $Y_W( z)$ a $\CC[[h]]$-module map
\begin{align*}
Y_W(z)\colon V\ot W&\to W((z))_h\\
v\ot w&\mapsto Y_W(z)(v\ot w)=Y_W(v,z)w=\sum_{r\in\mathbb{Z}} v_r w\ts z^{-r-1}
\end{align*}
which satisfies $Y_W(\vac,z)w=w$ for all $ w\in W$ and the weak associativity:
for any $u,v\in V$, $w\in W$ and $k \in\mathbb{Z}_{\geqslant 0}$
there exists $r\in\mathbb{Z}_{\geqslant 0}$
such that
\begin{align}
&\left((z_0 +z_2)^r\ts Y_W (u,z_0 +z_2)Y_W(v,z_2)\ts w \right.\non\\
&\left.\quad - (z_0 +z_2)^r\ts Y_W\big(Y(u,z_0)v,z_2\big)\ts w\right) 
\in  h^k   W[[z_0^{\pm 1},z_2^{\pm 1}]].\label{modulesassociativity}
\end{align}
\end{defn}
Regarding the definition, it is worth to recall that the weak associativity \eqref{modulesassociativity} implies
the  $\Sc$-Jacobi identity for the quantum vertex algebra module map $Y_W(\cdot ,z)$,
\begin{align}
&z_0^{-1}\delta\left(\frac{z_1 -z_2}{z_0}\right) Y_W(z_1)(1\ot Y_W(z_2))(u\ot v\ot w)\non\\
&\qquad-z_0^{-1}\delta\left(\frac{z_2-z_1}{-z_0}\right) Y_W(z_2)(1\ot Y_W(z_1)) \left(\Sc(-z_0)(v\ot u)\ot w\right)\non\\
\qquad\qquad=&\ts z_2^{-1}\delta\left(\frac{z_1 -z_0}{z_2}\right)Y_W(Y(u,z_0)v,z_2 )w\quad\text{for all}\quad u,v \in V\text{ and }w\in W.\label{wjacobi}
\end{align}
Indeed, the $\Sc$-Jacobi identity holds for the vertex operator map $Y(\cdot ,z)$, as it     is equivalent to weak associativity \eqref{associativity} and $\Sc$-locality \eqref{locality}; see \cite[Rem. 2.16]{Li}. Hence, one  deduces \eqref{wjacobi} from the $\Sc$-Jacobi identity for $Y(\cdot ,z)$ using  \cite[Lemma 5.7]{Li0} and \cite[Prop. 2.24]{Li}.

\begin{kor}\label{krlarkl}
Let $W$ be a restricted $\Drp$-module.
There exists a unique structure 
of $\vr$-module on $W$ such that 
\beq\label{rhstrana}
Y_W(L_{[n]}^+(u),z)=L_{[n]}(z+u)=\left(\R_{[n]}(v_1,\ldots ,v_n)L_1(v_1)\ldots L_n(v_n)\right)\Big|_{v_i=z+u_i}\Big. .
\eeq
\end{kor}

\begin{prf}
If $W$ is a restricted $\Drp$-module,  then \eqref{ghwqcd5} implies that the substitutions in \eqref{rhstrana} are well-defined and, furthermore, that the image of the module map $Y_W(\cdot, z)$ belongs to $\om(W,W((z))_h)$.
 One can now repeat the arguments from the corresponding part of the proof of Theorem \ref{mainn} to show that  \eqref{rhstrana}  defines a $\CC[[h]]$-module map on $\vr$  which satisfies the axioms from  Definition \ref{qvoamodule}.
\end{prf}

Roughly speaking,   the proof of Corollary \ref{krlarkl} is comprised of   showing that  a topologically free $\CC[[h]]$-module $W$ equipped with the map $L(u)\in \ndo\CC^N\ot\om(W,W((u))_h)$  satisfying the $RLL$-relation \eqref{rtt2} gives rise to a   $\vr$-module.
 The next proposition   establishes the converse.

\begin{pro}\label{jdnjdn}
Let $(W,Y_W)$ be a $\vr$-module. Then the operator series
\beq\label{sbsa2}
L(z)=Y_W(L^+(0)\vac,z) \in\ndo\CC^N\ot\om (W,W((z))_h)
\eeq
satisfies the $RLL$-relation  \eqref{rtt2}.
\end{pro}

\begin{prf}
The proposition is a simple consequence of the $\Sc$-Jacobi identity   for a $\vr$-module $W$. Let us apply the left-hand side of \eqref{wjacobi}  to 
\beq\label{utttz}
L_{23}^+(u)L_{14}^+(v)\R_{12}(-z_0+v-u)\ot w\quad\text{with }w\in W.
\eeq
 By using the explicit expression \eqref{Smap}  for the braiding,  we obtain
\begin{align*}
&z_0^{-1}\delta\left(\frac{z_1 -z_2}{z_0}\right)
Y_W(L_2^+(u),z_1)Y_W(L_1^+(v),z_2)\R_{12}(-z_0+v-u)w\non\\
& \qquad-z_0^{-1}\delta\left(\frac{z_2-z_1}{-z_0}\right) \R_{12}(-z_0+v-u) Y_W(L_1^+(v),z_2)  Y_W(L_2^+(u),z_1) w.
\end{align*}
Hence, by taking the residue with respect to $z_0$ we get
\begin{align}
& 
Y_W(L_2^+(u),z_1)Y_W(L_1^+(v),z_2)\R_{12}(-z_1+z_2+v-u)w\non\\
& \qquad- \R_{12}(z_2-z_1+v-u) Y_W(L_1^+(v),z_2)  Y_W(L_2^+(u),z_1) w.\label{utre4}
\end{align}

On the other hand, by using   formula \eqref{Ymap} for the vertex operator map, then \eqref{lemma21l} and, finally, the defining relations \eqref{rtt} we find
\begin{align*}
&Y(L_2^+(u),z_0)\ts L^+_1(v)\ts \R_{12}(-z_0+v-u)  
=  L_2 (z_0+u)\ts  L^+_1(v)\ts \R_{12}(-z_0+v-u)\\
=&\, \R_{12}(z_0+u-v)^{-1}L_2^+(z_0+u)\ts  L^+_1(v)\ts  \R_{12}(-z_0+v-u)=   L^+_1(v)\ts L_2^+(z_0+u).
\end{align*}
Therefore, by applying the right-hand side of the $\Sc$-Jacobi identity    to \eqref{utttz} we obtain
\beq\label{thisexpresom}
z_2^{-1}\delta\left(\frac{z_1 -z_0}{z_2}\right)Y_W(L^+_1(v)\ts L_2^+(z_0+u),z_2 )w.
\eeq
The residue of  \eqref{thisexpresom} with respect to the variable $z_0$  is zero. Thus, the two expressions in  \eqref{utre4} coincide. As $w\in W$ was arbitrary, by setting $u=v=0$ therein, we obtain the equality of operators on $W$,
$$
 \R_{12}(z_2-z_1) Y_W(L_1^+(0),z_2)  Y_W(L_2^+(0),z_1)   =Y_W(L_2^+(0),z_1)Y_W(L_1^+(0),z_2)\R_{12}(-z_1+z_2),
$$
as required.
\end{prf}

\section{Fixed points of the braiding}\label{sec07}

In this section, we study the braiding of $\vr$, as given by \eqref{Smap}.
Our main tools are two special cases of  the   fusion procedure  for the Yang $R$-matrix \eqref{betheyang} originated    in the work of A. Jucys \cite{J}; see also \cite[Sect. 6.4]{M} for more details.
Consider the usual action of the symmetric group $\mathfrak{S}_n$ on  $(\CC^N)^{\ot n}$, where   its elements permute the tensor factors.
 For any $n=1,\ldots ,N$ let $\Ec^+_{(n)}=H_{(n)}$ and $\Ec^-_{(n)}=A_{(n)}$  be the actions of the normalized symmetrizer and anti-symmetrizer,
\beq\label{symant}
h_{(n)}=\frac{1}{n!}\sum_{p\in \mathfrak{S}_n} p\in\CC[\mathfrak{S}_n]\Fand
 a_{(n)}=\frac{1}{n!}\sum_{p\in \mathfrak{S}_n} \sgn p \cdot p\in\CC[\mathfrak{S}_n],
\eeq
  on the tensor product space.
By \cite{J}, the consecutive evaluations $u_i = \pm (i-1)h $ for $i=1,\ldots ,n$ 
of the expression
$R_{[n]}(u_1,\ldots ,u_n)$ 
are well-defined and, furthermore, we have
\beq\label{bethefusion}
R_{[n]} (u_1,\ldots ,u_n)\big|_{u_1=0 }\big|_{u_2=\pm h }
 \dots \big|_{u_n=\pm (n-1)h }= n! \ts 
  \Ec^\pm_{(n)} .
\eeq
Applying this result to the normalized $R$-matrix $\R(u)=R^{\pm}(u)$ defined by \eqref{yang}, we get
\beq\label{bethefusionbethefusion}
R^\pm_{[n]} (u_1,\ldots ,u_n)\big|_{u_1=0 }\big|_{u_2=\pm h }
 \dots \big|_{u_n=\pm (n-1)h }=\alpha_n \ts 
  \Ec^\pm_{(n)} \quad\text{for}\quad
	\alpha_n= n!\prod_{1\leqslant i<j\leqslant n}\frac{j-i}{j-i+1}.
\eeq

For any $n=1,\ldots ,N$ consider the  series $L_{(n)}^\pm  (u)\in  V(R^{\pm}) [[u]]$ given by
\beq\label{obrisi}
L_{(n)}^\pm  (u)=\sum_{r\geqslant 1}l_{(n)}^\pm  (-r)u^{r-1}=\tr_{1,\ldots ,n}\ts   L_1^+( u  ) \ts L_2^+( u \pm h )\ldots L_n^+( u\pm (n-1)h ) \vac  .
\eeq

\begin{pro}\label{fiks}
The braiding \eqref{Smap} for $\vr=V(R^\pm)$ satisfies
\beq\label{idlt}
\Sc(z)(L^\pm_{ (n)} (u )\ot L^\pm_{ (m)} (v))
=L^\pm_{ (n)} (u )\ot L^\pm_{ (m)} (v)\quad\text{for all } m,n=1,\ldots ,N.
\eeq
\end{pro}

\begin{prf}
Let   $  w=(w_1,\ldots,w_n)$ and $ t= (t_1,\ldots ,t_m)$ be families of variables. We shall write 
$$
\wvr{u}=(u,u\pm h,\ldots ,u\pm (n-1)h)\fand \wvr{v}=(v,v\pm h,\ldots ,v\pm (m-1)h).
$$
By the defining relations  \eqref{rtt} for  $ \arpm$, we have
$$
L_{[n]}^+(w)=R^\pm_{[n]}(-w)^{-1}\ts \overleftarrow{L}_{[n]}^+(w)\ts R^\pm_{[n]}(w),\quad\text{where}\quad \overleftarrow{L}_{[n]}^+(w)=L_n^+(w_n)\ldots L_1^+(w_1).
$$
On the other hand, using the Yang--Baxter equation \eqref{ybe} and the unitarity  \eqref{uni}, one easily verifies the identity  $R^\pm_{[n]}(-w)^{-1}=R^\pm_{[n]}(w) $, so that the above equality turns to
$$
L_{[n]}^+(w)=R^\pm_{[n]}( w) \ts \overleftarrow{L}_{[n]}^+(w)\ts R^\pm_{[n]}(w).
$$
Applying the   evaluations $w_i=u\pm (i-1)h$  and using the fusion procedure \eqref{bethefusionbethefusion} we find
\beq\label{wfnd}
L_{[n]}^+(\wvr{u})=\alpha_n^2\ts \Ec^\pm_{(n)} \ts \overleftarrow{L}_{[n]}^+(\wvr{u})\ts \Ec^\pm_{(n)}. 
\eeq

The discussion above, along with explicit formula \eqref{Smap} for the action of braiding $\Sc(z)$, implies that  the left-hand  side of \eqref{idlt} equals
\begin{align}
&\tr_{1,\ldots, n +m}\ts
\ts   R_{n  m}^{\pm\ts \bar{1}\bar{2}}\ts L_{[n  ]}^{+13}(\bar{u}) L_{[m ]}^{+24}(\bar{v})\ts (R_{n  m }^{\pm\ts12})^{-1} \non\\
=\,&\tr_{1,\ldots, n +m}\ts
\ts   \alpha_n^2 \alpha_m^2\ts R_{n  m}^{\pm\ts\bar{1}\bar{2}}\ts \Ec^{\pm\ts 1}_{(n)} \ts   \Ec^{\pm\ts 2}_{(m)} \ts \overleftarrow{L}_{[n  ]}^{+13}(\bar{u}) \overleftarrow{L}_{[m ]}^{+24}(\bar{v})\ts \Ec^{\pm\ts 1}_{(n)} \ts   \Ec^{\pm\ts 2}_{(m)}\ts (R_{n  m }^{\pm\ts 12})^{-1},
\label{jtr56}
\end{align}
where  the idempotent   $\Ec^{\pm\ts 1}_{(n)} =\Ec^{\pm }_{(n)} $ (resp. $\Ec^{\pm\ts 2}_{(m)} =\Ec^{\pm\ts  }_{(m)}$) is applied on the tensor factors $1,\ldots ,n $  (resp. $n  +1,\ldots,n + m   )$ and we use the   notation  
\begin{gather*} 
R_{n  m}^{\pm\ts\bar{1}\bar{2}}=R_{n  m}^{\pm\ts\bar{1}\bar{2}}(z+\bar{u}-\bar{v})\fand
R_{n  m}^{\pm\ts 12}=R_{n  m }^{\pm\ts 12}(z+\bar{u}-\bar{v}).
\end{gather*}
By using the Yang--Baxter equation \eqref{ybe},   one easily verifies the equalities
\beq\label{tzh4}
R_{n  m}^{\pm\ts \bar{1}\bar{2}}\ts \Ec^{\pm\ts 1}_{(n)} \ts   \Ec^{\pm\ts 2}_{(m)}= \Ec^{\pm\ts 1}_{(n)} \ts   \Ec^{\pm\ts 2}_{(m)}\ts R_{n  m}^{\pm\ts 12}\fand
\Ec^{\pm\ts 1}_{(n)} \ts   \Ec^{\pm\ts 2}_{(m)}\ts (R_{n  m }^{\pm\ts 12})^{-1}=
 (R_{n  m }^{\pm\ts \bar{1}\bar{2}})^{-1}\ts \Ec^{\pm\ts 1}_{(n)} \ts   \Ec^{\pm\ts 2}_{(m)}.
\eeq
Hence, as $\Ec^{\pm }_{(n)}$ and $   \Ec^{\pm }_{(m)}$ are idempotents,  the expression in \eqref{jtr56} coincides with
$$
\tr_{1,\ldots, n +m}\ts
\ts   \alpha_n^2 \alpha_m^2\ts \Ec^{\pm\ts 1}_{(n)} \ts   \Ec^{\pm\ts 2}_{(m)}\ts    R_{n  m}^{\pm\ts 12}\ts \Ec^{\pm\ts 1}_{(n)} \ts   \Ec^{\pm\ts 2}_{(m)}\ts \overleftarrow{L}_{[n  ]}^{+13}(\bar{u}) \overleftarrow{L}_{[m ]}^{+24}(\bar{v})\ts (R_{n  m }^{\pm\ts \bar{1}\bar{2}})^{-1}\ts \Ec^{\pm\ts 1}_{(n)} \ts   \Ec^{\pm\ts 2}_{(m)}.
$$
Next, we employ the cyclic property of the trace to move the leftmost copies of $\Ec^{\pm }_{(n)}$ and $   \Ec^{\pm }_{(m)}$ to the right, thus getting
$$
\tr_{1,\ldots, n +m}\ts
\ts   \alpha_n^2 \alpha_m^2\ts     R_{n  m}^{\pm\ts 12}\ts \Ec^{\pm\ts 1}_{(n)} \ts   \Ec^{\pm\ts 2}_{(m)}\ts \overleftarrow{L}_{[n  ]}^{+13}(\bar{u}) \overleftarrow{L}_{[m ]}^{+24}(\bar{v})\ts (R_{n  m }^{\pm\ts \bar{1}\bar{2}})^{-1}\ts \Ec^{\pm\ts 1}_{(n)} \ts   \Ec^{\pm\ts 2}_{(m)}.
$$
Once again,  we make use of the second identity in \eqref{tzh4} and move the rightmost copies of $\Ec^{\pm }_{(n)}$ and $   \Ec^{\pm }_{(m)}$ to the left, so that we obtain
$$
\tr_{1,\ldots, n +m}\ts
\ts  \alpha_n^2 \alpha_m^2\ts     R_{n  m}^{\pm\ts 12}\ts \Ec^{\pm\ts 1}_{(n)} \ts   \Ec^{\pm\ts 2}_{(m)}\ts \overleftarrow{L}_{[n  ]}^{+13}(\bar{u}) \overleftarrow{L}_{[m ]}^{+24}(\bar{v})\ts \Ec^{\pm\ts 1}_{(n)} \ts   \Ec^{\pm\ts 2}_{(m)}\ts (R_{n  m }^{\pm\ts 12})^{-1}.
$$
Finally, this is equal to
$$
\tr_{1,\ldots, n +m}\ts
\ts   \alpha_n^2 \alpha_m^2\ts   \Ec^{\pm\ts 1}_{(n)} \ts   \Ec^{\pm\ts 2}_{(m)}\ts \overleftarrow{L}_{[n  ]}^{+13}(\bar{u}) \overleftarrow{L}_{[m ]}^{+24}(\bar{v})\ts \Ec^{\pm\ts 1}_{(n)} \ts   \Ec^{\pm\ts 2}_{(m)} 
$$
as, due to the cyclic property of the trace, the terms $R_{n  m}^{\pm\ts 12}$ and $(R_{n  m }^{\pm\ts 12})^{-1}$  cancel.  It remains to observe that  by \eqref{wfnd}, this coincides with
the right-hand side of \eqref{idlt}, as required.
\end{prf}

Let $\pi^\pm\colon A(R^\pm)\to \Pc^\pm$ be the map   $a\mapsto a\cdot 1$, where the action of $a\in A(R^\pm)$ on $1\in\Pc^\pm$ is given by Proposition \ref{repje}. For any $n=1,\ldots ,N$ and $r\geqslant 1$ introduce the polynomials
$$
t(n,r)=\sum_{\substack{r_1,\ldots,r_n\geqslant 1\\r_1+\ldots+r_n=n+r-1 }}
t_{r_1} t_{r_2}\ldots t_{r_n}\in\CC[t_1,\ldots ,t_n].
$$

\begin{pro}
For any $n=1,\ldots ,N$ and $r=1,2,\ldots$ we have
\begin{align*}
&\pi^+(l_{(n)}^+ (-r))=  t(n,r)\sum_{  i_1,\ldots,i_n=1,\ldots,N} x_{i_1} y_{i_1}\ldots x_{i_n} y_{i_n}\mod h,\\
&\pi^-(l_{(n)}^- (-r))=n!\ts\ts t(n,r)\sum_{1\leqslant i_1<\ldots <i_n\leqslant N} x_{i_1} y_{i_1}\ldots x_{i_n} y_{i_n}\mod h.
\end{align*}
Furthermore, for any $n=1,\ldots ,N$ the coefficients $l_{(n)}^+ (-r)$ with $r=1,2,\ldots$ are algebraically independent.
\end{pro}

\begin{prf}
The images of the coefficients of the series $L_{(n)}^\pm (u)$ under the map $\pi^\pm$   are easily computed modulo $h$ using formula \eqref{obrisi}. As for the second assertion, 
the algebraic independence of any $l_{(n)}^+ (-s_1),\ldots ,l_{(n)}^+ (-s_m)$
with $   s_1>\ldots > s_m\geqslant 1$ is established by considering the partial derivatives
$$
\frac{\partial }{\partial t_{s_j  }}\pi^+(l_{(n)}^+ (-s_i))\mod h \qquad\text{for }i,j=1\ldots ,m
$$
and using the Jacobian Criterion.
\end{prf}

 The next theorem relies on  the   property of the Yang $R$-matrix,
\beq\label{bethefsi}
A_{(N)}\ts R_{0N}(u+(N-1)h)\ldots R_{02}(u+h)\ts R_{01}(u) =A_{(N)}\left(1-\frac{h}{u}\right); 
\eeq
 see \cite[Ch. 1]{M} for more details.
It strengthens Proposition \ref{fiks}   in the case of the series
$$\wvr{\mathbb{L}} (u)\coloneqq L_{(N)}^{-}  (u)  \in  V(R^{-}) [[u]].$$

\begin{thm}\label{ratfiks} 
For any $a\in V(R^{-})$ we have
\beq\label{idlt6}
\Sc(z)(a\ot \wvr{\mathbb{L}} (v))
=a\ot\wvr{\mathbb{L}}  (v).
\eeq
\end{thm}

\begin{prf}
Clearly, it is sufficient to check that \eqref{idlt6} holds for $a=L_{[n]}^+(u)$ with $n\geqslant 1$ and $u=(u_1,\ldots ,u_n)$ as the case $a=\vac$ follows from the general theory; see \cite[Prop. 1.6]{EK}.
Let us write 
$ 
 \wvr{v}=(v , v-h\ldots ,v- (N-1) h).
$
By the definition \eqref{Smap} of   braiding,      the left-hand side of \eqref{idlt6} with $a=L_{[n]}^+(u)$ equals
\begin{align*}
\Sc(z)(L_{[n]}^+(u)\ot \wvr{\mathbb{L}} (v))=
 \tr_{n+1,\ldots, n +N} 
\ts   R_{n  N}^{-\ts \bar{1}\bar{2}}(z+u-\bar{v})\ts L_{[n  ]}^{+13}(u) L_{[N ]}^{+24}(\bar{v})\ts R_{n  N }^{-\ts 12}(z+u-\bar{v})^{-1}. 
\end{align*}
Arguing as in the proof of Proposition \ref{fiks} (cf. \eqref{wfnd}), one shows that this equals
\begin{align}
 \tr_{n+1,\ldots, n +N}
\ts    R_{n  N}^{-\ts \bar{1}\bar{2}}(z+u-\bar{v})\ts  A_{(N)}\ts L_{[n  ]}^{+13}(u) L_{[N ]}^{+24}(\bar{v})\ts R_{n  N }^{-\ts 12}(z+u-\bar{v})^{-1}, \label{hjkl7}
\end{align}
where  the anti-symmetrizer $A_{(N)}$ is applied on the tensor factors $n+1,\ldots ,n+N $. 
Using the Yang--Baxter equation \eqref{ybe}   one easily proves the equality
\beq\label{itsrfsn}
R_{n  N}^{-\ts \bar{1}\bar{2}}(z+u-\bar{v})\ts  A_{(N)}=
   A_{(N)} \ts R_{n  N}^{-\ts \bar{1}2}(z+u-\bar{v}) .
\eeq
Next, by   the property \eqref{bethefsi} of the anti-symmetrizer,   the right-hand side of \eqref{itsrfsn} equals $A_{(N)}F(z,u,v)$, where the power series $ F(z,u,v)$ is given by
$$
F(z,u,v)= \prod_{i=1}^n \left(1-\frac{h}{z+u_i-v}\right) \prod_{j=1}^N \left(1+\frac{h}{z+u_i-v+(j-1)h}\right)^{-1}. 
$$
Hence, \eqref{hjkl7} coincides with
\beq\label{thtrm6}
\tr_{n+1,\ldots, n +N}
\ts    F (z,u,v)  A_{(N)}\ts L_{[n  ]}^{+13}(u) L_{[N ]}^{+24}(\bar{v})\ts R_{n  N }^{-\ts 12}(z+u-\bar{v})^{-1}. 
\eeq

We now slightly modify this argument to remove the term $R_{n  N }^{-\ts 12}(z+u-\bar{v})^{-1}$.
First,     in analogy with  the first part of the proof of   Proposition \ref{fiks} (cf. \eqref{wfnd}), we have  
\beq\label{gdawr1}
A_{(N)}\ts  L_{[N ]}^{+24}(\bar{v})
=\alpha_N^2\ts A_{(N)}\ts A_{(N)}\ts  \cev{L}_{[N ]}^{+24}(\bar{v})\ts A_{(N)} 
=\alpha_N^2\ts A_{(N)}\ts  \cev{L}_{[N ]}^{+24}(\bar{v})\ts A_{(N)} .
\eeq
Next, we use  the property \eqref{bethefsi} of the anti-symmetrizer once again to obtain
\beq\label{gdawr2}
A_{(N)}\ts R_{n  N }^{-\ts 12}(z+u-\bar{v})^{-1}= A_{(N)}\ts F(z,u,v)^{-1}.
\eeq
Finally, note that, by the fusion procedure \eqref{bethefusionbethefusion}, the right-hand side of \eqref{gdawr1} equals  $L_{[N ]}^{+24}(\bar{v})$.
Hence,
we can use \eqref{gdawr1} and  \eqref{gdawr2} to rewrite \eqref{thtrm6} as
$$
\tr_{n+1,\ldots, n +N}
\ts      L_{[n  ]}^{+13}(u) L_{[N ]}^{+24}(\bar{v}) 
= L_{[n]}^+(u)\ot\wvr{\mathbb{L}}(v),
$$
which concludes the proof.
\end{prf}

\section{Commutative families   in the Yangian quantization of   $\mathcal{O}(\mathfrak{gl}_N((z^{-1})))$}\label{section06}


In this section, we study commutative families and central elements in a certain algebra, which are closely related with the families of fixed points of the braiding established in the previous section. The aforementioned algebra  was introduced by Krylov and Rybnikov  and it can be regarded as 
the Yangian quantization of the Poisson algebra $\mathcal{O}(\mathfrak{gl}_N((z^{-1})))$; see \cite{KR} for more details.
Consider the  following normalization of the Yang $R$-matrix \eqref{betheyang}:
\beq\label{chsck}
\wht{R}(u) =u R(u)=uI-hP\in\ndo\CC^N\ot \ndo\CC^N[ u,h ].
\eeq
Its inverse is given by
\beq\label{betheinv}
\wht{R}(u)^{-1}=\frac{\wht{R}(-u)}{u^2-h^2}  = (uI+hP)\sum_{l\geqslant 0}\frac{h^{2l}}{u^{2l+2}} \in\ndo\CC^N\ot \ndo\CC^N[ u^{-1}][[h] ].
\eeq

We follow \cite[Sect. 3.1]{KR} to introduce the   algebra $\Yht$ and its completion.  The algebra is defined over the commutative ring $\CC[[h]]$ and it is obtained from the original definition, which is given over the complex field, by suitably rescaling its generators and the spectral parameter.  Let $\Yht$ be the   associative algebra over  $\CC[[h]]$  generated by the elements $\lambda_{ij}^{(r)}$, where $r\in\ZZ$  and $i,j=1,\ldots ,N$,   subject to the defining relations 
\beq\label{bethertt}
\wht{R} ( u-v )\ts  \Lc_1(u)\ts \Lc_2 (v)
=  \Lc_2 (v)\ts \Lc_1(u)\ts\wht{R} ( u-v  ).
\eeq
The matrix of generators $\Lc (u) $ is defined by
$$
\Lc (u) =\sum_{i,j=1}^N e_{ij}\ot \lambda_{ij}  (u),\quad\text{where}\quad \lambda_{ij} (u)= \sum_{r\in\ZZ} \lambda_{ij}^{(-r)}u^{r }.
$$

We now establish a connection between $\Yht$-modules and  $\vr$-modules. Consider the operator series $L(u)=Y(L^+(0),u)$; recall \eqref{preslikavanjeel} and \eqref{Ymap}. It is evident from \eqref{rtt2} that it satisfies the $RLL$-relation \eqref{bethertt}. Indeed, the relation is    obtained by  multiplying  \eqref{rtt2} by the polynomial $u-v\mp h$. 
Therefore, by Proposition \ref{jdnjdn} we have

\begin{kor}\label{rktrhtr}
Let $(W,Y_W)$ be a $\vr$-module. Then the assignment
$$
\Lc(z)\mapsto L(z)=Y_W(L^+(0)\vac,z)
$$
defines a structure of $\Yht$-module on $W$ such that
$$
\Lc(z) \in\ndo\CC^N\ot\om (W,W((z))_h).
$$
In particular, the above assignment defines an action of the algebra $\Yht$ on $\vr$.
\end{kor}

Next, we complete the algebra $\Yht$; cf. \cite[Sect. 3.1]{KR}. Let $I_p$ with $p\geqslant 0$ be the two-sided ideal in  $\Yht$ generated by $h^{p+1}$ and all $\lambda_{ij}^{(-r)}$ with $r>p$. Define the completed algebra  by
$$
\Yhtc =\lim_{\longleftarrow} \Yht/I_p.
$$

Multiplying    relation \eqref{bethertt} by the inverse $\wht{R} ( u-v )^{-1}$ from the left  we find
\beq\label{bethecomp}
   \Lc_1(u)\ts \Lc_2 (v)
=  \wht{R} ( u-v )^{-1}  \Lc_2 (v)\ts \Lc_1(u)\ts\wht{R} ( u-v  ).
\eeq
Indeed,  \eqref{betheinv} implies 
$$\wht{R} ( u-v )^{-1}\in 
\ndo\CC^N \ot \ndo\CC^N [v][[u^{-1},h]]
,$$ 
while, for any $p\geqslant 0$, the expression $\Lc_1(u) \Lc_2 (v) $, when regarded modulo $I_p$, belongs to 
$$ 
\ndo\CC^N \ot \ndo\CC^N \ot \Yht ((u^{-1},v^{-1})) .
$$ 
 Thus, the aforementioned product is well-defined so that it  gives rise to   \eqref{bethecomp}.

For any   $n=1,\ldots ,N$ consider the series 
\beq\label{bethegaqw}
\Lc_{(n)}^\pm (z)=
\tr_{1,\ldots ,n}\ts \Ec_{(n)}^\pm\ts\Lc_1(z )\Lc_2(z\pm h)\ldots \Lc_n(z\pm (n-1) h)\in \Yhtc[[z^{\pm 1}]].
\eeq
The next proposition was proved in \cite[Prop. 3.15]{KR}  in the case of the  anti-symmetrizer $\Ec_{(n)}^- =A_{(n)}$. Although the proof of  its generalization to the case of the symmetrizer $\Ec_{(n)}^+ =H_{(n)}$  relies on analogous techniques, we present its details in both cases so that  we can refer to them     later on.

\begin{pro}\label{bethemc}
The coefficients of  all $\Lc_{(n)}^\pm (z)$
  mutually commute.
\end{pro}

\begin{prf}
(1) First, we establish some   preliminary results which we shall need in the main part of the proof.
 Choose any $n_1,n_2=1,\ldots ,N$. Denote by $\Ec_{(n_i)}$ for $i=1,2$ the symmetrizer $H_{(n_i)}$ or the anti-symmetrizer $A_{(n_i)}$. 
For $i=1,2$ introduce the families of variables
$$ 
x_i=(z_i+c_1^{(i)}h,z_i + c_2^{(i)}h\ldots ,z_i +c_{n_i}^{(i)}h),\quad\text{where}\quad 
c_a^{(i)}=\begin{cases}
a-1,&\text{if }\Ec_{(n_i)}=H_{(n_i)}\\
1-a,&\text{if }\Ec_{(n_i)}=A_{(n_i)} 
\end{cases}
$$
and $z_1,z_2$ are single variables.
The Yang--Baxter equation \eqref{ybe}  and   \eqref{bethefusion} imply   
\begin{align}
\Ec_{(n_1)}^1\ts \Ec_{(n_2)}^2\ts  \wht{R}_{n_1 n_2}^{ 12}( x_1-x_2    )
&=\wht{R}_{n_1 n_2}^{\bar{\scriptstyle 1}\bar{\scriptstyle 2}}(  x_1-x_2    ) \ts
\Ec_{(n_1)}^1\ts \Ec_{(n_2)}^2 ,\label{betheeyb2}
\end{align}
(cf. \eqref{tzh4}), where    $\Ec^1_{(n_1)}$ (resp. $\Ec^2_{(n_2)}$) indicates that the idempotent $ \Ec_{(n_1)}$ (resp. $ \Ec_{(n_2)}$) is applied on the tensor factors $1,\ldots ,n_1$ (resp. $n_1+1,\ldots ,n_1+n_2$).
Next, by using the Yang--Baxter equation and    \eqref{bethecomp}   one can verify the identity
\beq\label{thrltn}
   \Lc_{[n]}^{13}(u)\ts \Lc_{[m]}^{23} (v)
=
\wht{R}_{nm}^{\bar{\scriptstyle 1}\bar{\scriptstyle 2}}( u-v )^{-1}
\Lc_{[m]}^{23} (v)\ts \Lc_{[n]}^{13}(u)\ts \wht{R}_{nm}^{ 12}(  u-v   ), 
\eeq
where $u=(u_1,\ldots ,u_n)$ and $v=(v_1,\ldots ,v_m)$ are families of variables and
\beq\label{bethetzh4}
\Lc_{[n]}(u)=\frac{1}{n!}R_{[n]}(u_1,\ldots ,u_n)\ts  \Lc_{1}(u_1)\ldots \Lc_{n}(u_n).
\eeq
Finally, by the fusion procedure \eqref{bethefusion}, applying the consecutive evaluations $u_a=z_i+c_i^{(a)} h$  and then the trace $\tr_{1,\ldots ,n_i}$ to \eqref{bethetzh4} with $n=n_i$, we obtain   $\Lc_{(n_i)}  (z_i)\coloneqq \Lc_{(n_i)}^\pm (z_i)$, where, as before, the plus (resp. minus) sign corresponds to the symmetrizer (resp. anti-symmetrizer) case.\vspace{2pt}

\noindent (2) We are now prepared for the main part of the proof.
We shall write $k=n_1+n_2$, 
\beq\label{bethenott} 
\wndr{\Lc}_{[n]}(u)=\Lc_{1}(u_1)\ldots \Lc_{n}(u_n)\Fand \overleftarrow{\wndr{\Lc}}_{[n]}(u)=\Lc_{n}(u_n)\ldots \Lc_{1}(u_1).
\eeq
Clearly, it is sufficient to check that  
$ 
\Lc_{(n_1)}(z_1)\Lc_{(n_2)}(z_2)
$    coincides with $ 
\Lc_{(n_2)}(z_2)\Lc_{(n_1)}(z_1)
$.
By \eqref{bethegaqw}, the former    is equal to
\begin{align*}
\tr_{1,\ldots ,k} \ts \Ec^1_{(n_1)} \ts \Lc_1(z_1+c^{(1)}_1h) \ldots \Lc_{n_1}(z_1+c^{(1)}_{n_1}h)
\ts \Ec^2_{(n_2)} \ts \Lc_{n_1+1}(z_2+c^{(2)}_1h) \ldots \Lc_{k}(z_2 +c^{(2)}_{n_2}h).
\end{align*}  
Using the fusion procedure \eqref{bethefusion}, we rewrite this expression as
$$
 \tr_{1,\ldots, k} \ts    \Lc^{13}_{[n_1]}(x_1)  
\ts \Lc^{23}_{[n_2]}(x_2).
$$
The $RLL$-relation \eqref{thrltn} implies that it is equal to 
\begin{align}
&\, \tr_{1,\ldots,  k}   
\ts  \wht{R}_{n_1 n_2}^{\bar{\scriptstyle 1}\bar{\scriptstyle 2}}( x_1-x_2  )^{-1}\Lc^{23}_{[n_2]}(x_2)\ts   \Lc^{13}_{[n_1]}(x_1)\ts \wht{R}_{n_1 n_2}^{12}( x_1-x_2  )\non\\
=&\,  \tr_{1,\ldots , k}   
\ts \wht{R}_{n_1 n_2}^{\bar{\scriptstyle 1}\bar{\scriptstyle 2}}( x_1-x_2  )^{-1}\ts
 \Ec^2_{(n_2)}\ts
\wndr{\Lc}^{23}_{[n_2]}(x_2)\ts \Ec^1_{(n_1)}\ts  \wndr{\Lc}^{13}_{[n_1]}(x_1)\ts \wht{R}_{n_1 n_2}^{12}( x_1-x_2  )\non\\ 
=&\,  \tr_{1,\ldots , k}   
\ts \wht{R}_{n_1 n_2}^{\bar{\scriptstyle 1}\bar{\scriptstyle 2}}(x_1-x_2  )^{-1}\ts
\Ec^1_{(n_1)}\ts \Ec^2_{(n_2)}\ts
\Ec^1_{(n_1)}\ts \Ec^2_{(n_2)}\ts
\wndr{\Lc}^{23}_{[n_2]}(x_2)\ts   \wndr{\Lc}^{13}_{[n_1]}(x_1)\ts \wht{R}_{n_1 n_2}^{12}( x_1-x_2  ),\label{bethetrfr1}
\end{align}
where, in the last equality, we employed the fact that $\Ec_{(n_i)}$ are idempotents.   
Let us use \eqref{betheeyb2}  to move the leftmost copies of $\Ec^1_{(n_1)}  \Ec^2_{(n_2)}$ in \eqref{bethetrfr1}  to the left. Then, we use \eqref{betheeyb2} and
\beq\label{betheeyb3}
\Ec^i_{n_i} \ts \wndr{\Lc}_{[n_i]}(x_i) =\overleftarrow{\wndr{\Lc}}_{[n_i]}(x_i)\ts \Ec^i_{n_i}\quad\text{for }i=1,2 
\eeq
to move the rightmost copies of $\Ec^1_{(n_1)}  \Ec^2_{(n_2)}$ to the right.
Thus,   \eqref{bethetrfr1}  becomes
$$
  \tr_{1,\ldots, k}  \ts
 \Ec^1_{(n_1)}  \ts\Ec^2_{(n_2)} 
\ts \wht{R}_{n_1 n_2}^{12}( x_1-x_2 )^{-1}
\overleftarrow{\wndr{\Lc}}^{23}_{[n_2]}(x_2)\ts   \overleftarrow{\wndr{\Lc}}^{13}_{[n_1]}(x_1)\ts \wht{R}_{n_1 n_2}^{\bar{\scriptstyle 1}\bar{\scriptstyle 2}}( x_1-x_2 )\ts
 \Ec^1_{(n_1)}  \ts\Ec^2_{(n_2)} .
$$
By using the cyclic property of the trace, we can move the leftmost copies of the idempotents $\Ec^1_{(n_1)}$ and $  \Ec^2_{(n_2)}$ to the right, thus getting
$$
  \tr_{1,\ldots, k}  \ts
   \wht{R}_{n_1 n_2}^{12}( x_1-x_2 )^{-1}
\overleftarrow{\wndr{\Lc}}^{23}_{[n_2]}(x_2)\ts   \overleftarrow{\wndr{\Lc}}^{13}_{[n_1]}(x_1)\ts \wht{R}_{n_1 n_2}^{\bar{\scriptstyle 1}\bar{\scriptstyle 2}}( x_1-x_2 )\ts
 \Ec^1_{(n_1)}  \ts\Ec^2_{(n_2)} .
$$
Moreover, we use \eqref{betheeyb2} and \eqref{betheeyb3} to   move the remaining copies back to the left:
$$
  \tr_{1,\ldots, k}  \ts
   \wht{R}_{n_1 n_2}^{12}( x_1-x_2 )^{-1} \Ec^2_{(n_2)} \ts
 \wndr{\Lc}^{23}_{[n_2]}(x_2)\ts    \Ec^1_{(n_1)}\ts  \wndr{\Lc}^{13}_{[n_1]}(x_1)\ts \wht{R}_{n_1 n_2}^{12}( x_1-x_2 ) 
  .
$$
Finally, by the cyclic property of the trace the $R$-matrices on the left and    the right cancel, so that   the   expression turns to
$$
  \tr_{1,\ldots, k}  \ts
     \Ec^2_{(n_2)} \ts
 \wndr{\Lc}^{23}_{[n_2]}(x_2)\ts    \Ec^1_{(n_1)}\ts  \wndr{\Lc}^{13}_{[n_1]}(x_1)     =
\Lc_{(n_2)}(z_2)\ts \Lc_{(n_1)}(z_1)
 ,
$$
as required.
\end{prf}

Let us compare Propositions \ref{fiks} and   \ref{bethemc}.
First of all, note that  the   analogue of the series $\Lc_{(n)}^\pm (z)$  for the quantum vertex algebra $\vr$ or its modules (in the sense of Corollary \ref{rktrhtr}) is equal, up to a nonzero multiplicative scalar from $\CC$, to
\beq\label{lbblebe}
\left(
\tr_{1,\ldots ,n}\ts \R_{[n]}(u_1,\ldots ,u_n)L_1(u_1)\ldots L_n(u_n)
\right)\Big|_{u_i=z\pm (i-1)h},
\eeq
where the substitutions are taken for all $i=1,\ldots ,n$.
On the other hand, \eqref{lbblebe} is exactly the image
$ 
Y(L_{(n)}^{\pm}(0),z)
$ 
(resp. $Y_W(L_{(n)}^{\pm}(0),z)$)
of the series \eqref{obrisi} at $u=0$ under the vertex operator map \eqref{Ymap} (resp. $\vr$-module map $Y_W(\cdot ,z)$). Finally, recall that by Proposition \ref{fiks} the coefficients of the aforementioned series give rise to   fixed points of the braiding. Hence, in view of the partial refinement of   Proposition \ref{fiks}, as given by Theorem \ref{ratfiks}, it is natural to investigate whether the same can be done with Proposition \ref{bethemc}. Such a result is given in the following theorem.
Consider the  special case of   series \eqref{bethegaqw},
\beq\label{bethert}
\mathbb{L}(z) \coloneqq \mathcal{L}_{(N)}^-(z)=
\tr_{1,\ldots ,N}\ts A_{(N)}\ts\Lc_1(z )\Lc_2(z- h)\ldots \Lc_N(z-(N-1) h) .
\eeq
 
\begin{thm}\label{centralnithm}
The coefficients of    $\mathbb{L}(z)$
belong to the center of $\Yhtc$.
\end{thm}

\begin{prf}
It is sufficient to check the equality $\mathcal{L}(z_0)\mathbb{L}(z)=\mathbb{L}(z)\mathcal{L}(z_0)$. Using \eqref{bethert} we rewrite its left-hand side as
$$
\mathcal{L}(z_0)\ts \mathbb{L}(z)=\tr_{1,\ldots ,N}\ts A_{(N)}\ts\mathcal{L}_0(z_0)\ts \Lc_1(z )\Lc_2(z- h)\ldots \Lc_N(z-(N-1) h) .
$$
Next, we employ \eqref{bethecomp} to move $\mathcal{L}_0(z_0)$ to the right, thus getting
\beq\label{betheztp}
\tr_{1,\ldots ,N}\ts  A_{(N)}\ts\wht{R}^{-1}\ts \Lc_1(z )\Lc_2(z- h)\ldots \Lc_N(z-(N-1) h) \ts \mathcal{L}_0(z_0)\ts\wht{R},
\eeq
where 
$$
\wht{R}=\wht{R}_{0N}(z_0-z+(N-1)h)\ldots \wht{R}_{02}(z_0-z+ h)\wht{R}_{01}(z_0-z ).
$$
By \eqref{bethefsi}, we have
\beq\label{bethesd}
A_{(N)}\ts\wht{R}^{\pm 1}= A_{(N)} \ts f(z_0-z)^{\pm 1},
\quad\text{where}\quad
f(x)=(x-h) \prod_{i= 1,2,\ldots,N-1}(x+ih)  .
\eeq
On the other hand, by combining the fusion procedure \eqref{bethefusion}   and   relation \eqref{bethertt} we find
\beq\label{bethesdd}
A_{(N)}\ts \Lc_1(z) \ldots \Lc_N(z-(N-1)h) =\Lc_N(z-(N-1)h)\ldots   \Lc_1(z) \ts A_{(N)}.
\eeq

Consider the expression in \eqref{betheztp}. One can use \eqref{bethesd} and \eqref{bethesdd} to move the anti-symmetrizer all the way to the right, thus canceling all $R$-matrix factors, and then to return it  back to the left, thus getting
$$
\tr_{1,\ldots ,N}\ts A_{(N)}\ts \Lc_1(z )\Lc_2(z- h)\ldots \Lc_N(z-(N-1) h) \ts\mathcal{L}_0(z_0)=\mathbb{L}(z)\ts \mathcal{L}(z_0),
$$
which completes the proof.
\end{prf}

\section{Trigonometric case}\label{zadnjisection}

In this section, we  discuss a generalization of results from Section \ref{section06} to the case of  the trigonometric $R$-matrix in type $A$. 
We  use the calligraphic font for the trigonometric $R$-matrix in order to easily distinguish it  from its rational counterpart  in the text. Introduce the two-parameter $R$-matrix $\Rc (x,y)\in\ndo\CC^N\ot\ndo\CC^N[[h]][x,y]$,
\begin{align}
\Rc (x,y) =&\left(xe^{-h/2}-ye^{h/2}\right)\sum_{i=1}^N e_{ii}\ot e_{ii} 
+ (x-y)\sum_{\substack{i,j=1\\i\neq j}}^N e_{ii}\ot e_{jj} \non\\
&+ \left(e^{-h/2}-e^{h/2}\right)x \sum_{\substack{i,j=1\\i> j}}^N e_{ij}\ot e_{ji}
+\left(e^{-h/2}-e^{h/2}\right)y\sum_{\substack{i,j=1\\i< j}}^N e_{ij}\ot e_{ji}. \label{twotrigR}
\end{align}
Next, let $P^h $ be  the $h$-permutation operator,
$$
P^h = \sum_{i=1}^N e_{ii}\ot e_{ii} + e^{h/2}\sum_{\substack{i,j=1\\i> j}}^N e_{ij}\ot e_{ji} +e^{-h/2}\sum_{\substack{i,j=1\\i< j}}^N e_{ij}\ot e_{ji} .
$$
Consider the action of the symmetric group $\mathfrak{S}_n$ on the space $(\CC^N)^{\ot n}$ which is uniquely determined by the requirement that the transpositions $(i,i+1)\in\mathfrak{S}_n $ act  as the $h$-permutation operator $P^h$ on the tensor factors $i$ and $i+1$, i.e. $(i,i+1)=P^h_{i\ts i+1}$.
Let $A^h_{(n)}$ be the image of the normalized anti-symmetrizer $a_{(n)}\in \CC[\mathfrak{S}_n]$, as given by \eqref{symant}, with respect to this action. Finally, we recall  the fusion procedure \cite{C} for  the $R$-matrix  \eqref{twotrigR}:
\beq\label{trigfusion}
\prod_{i=1,\ldots ,n-1}^{\longrightarrow}\prod_{j=i+1,\ldots ,n}^{\longrightarrow} \Rc_{ij}(xe^{-(i-1)h},xe^{-(j-1)h})
=
n!\, x^{\frac{n(n-1)}{2}}
\prod_{0\leqslant i< j\leqslant n-1} (e^{-ih}-e^{-jh})\ts A^h_{(n)}.
\eeq

We shall also need the one-parameter $R$-matrix  
\beq\label{rhttrg}
\wht{\Rc} (x )=e^{-h/2}\ts \Rc(x,1) \in \ndo\CC^N\ot\ndo\CC^N[[h]] [x] 
\eeq
It  satisfies the    Yang--Baxter equation  
\beq\label{trigybe}
\wht{\Rc}_{12}(x)\ts \wht{\Rc}_{13}(xy)\ts \wht{\Rc}_{23}(y)=\wht{\Rc}_{23}(y)\ts \wht{\Rc}_{13}(xy)\ts \wht{\Rc}_{12}(x) .
\eeq
Moreover, we have 
\beq\label{hattriguni}
 \wht{\Rc}_{12}(x)\ts  \wht{\Rc}_{21}(1/x) =\left(e^{-h}x-1\right)\left(e^{-h}x^{-1}-1\right),
\eeq
so that its inverse $\wht{\Rc}(x)^{-1}$ belongs  to $\ndo\CC^N\ot \ndo\CC^N[[x,h]]$. 

\begin{rem}\label{remarkrmatricca}
The $R$-matrix $\wht{\Rc} (x)$ is closely related with  the rational $R$-matrix $\wht{R} (u)$ defined by \eqref{chsck}.
More specifically, by setting $x=e^u$ in $\wht{\Rc}(x)$   we obtain the $R$-matrix   
$$\wht{\Rc}(e^u)= \wht{\Rc}(x)\big|_{x=e^u}\in\ndo\CC^N \ot\ndo\CC^N  [[u,h]] .$$  
Consider the $\ZZ$-gradation on $\ndo\CC^N \ot\ndo\CC^N [u ,h]$ given by $\deg u^r h^s =-r-s$. Extend the degree function to
$\ndo\CC^N \ot\ndo\CC^N  [[u,h]]$ by allowing it to take the infinite value. The $R$-matrix $\wht{\Rc}(e^u)$ is of finite degree and, furthermore, its degree is $-1$. Finally, its component of degree $-1$ is the rational $R$-matrix $\wht{R}(u)$. In addition, by rewriting the identity established in the proof of  \cite[Lemma 4.3]{FJMR} in terms of $\wht{\Rc} (x )$, one obtains the  trigonometric counterpart of the property \eqref{bethefsi} of the Yang $R$-matrix,
\beq\label{qantisimprop}
A^h_{(N)}\wht{\Rc}_{0N}(xe^{(N-1)h} )\ldots \wht{\Rc}_{02}(xe^{ h} )
\wht{\Rc}_{01}(x  ) 
=A^h_{(N)}
\frac{   xe^{-h }-1}{e^{ (N-1)h/2}}\prod_{i=2}^N \left(x e^{(i-1)h} -1\right).
\eeq
\end{rem}

Motivated by \cite[Sect. 3]{FPT}, we consider the associative algebra $\tYht$ over the commutative ring $\CC[[h]]$ which is generated by the elements $\tau_{ij}^{(r)}$, where $r\in\ZZ$  and $i,j=1,\ldots ,N$.   Its defining relations are given by
\beq\label{trigbethertt}
\wht{\Rc} ( x/y )\ts  \Tc_1(x)\ts \Tc_2 (y)
=  \Tc_2 (y)\ts \Tc_1(x)\ts\wht{\Rc} ( x/y  ),
\eeq
where the matrix of generators $\Tc (x) $ is given by
$$
\Tc (x) =\sum_{i,j=1}^N e_{ij}\ot \tau_{ij}  (x) \quad\text{with}\quad \tau_{ij} (x)= \sum_{r\in\ZZ} \tau_{ij}^{(-r)}x^{r }.
$$

In parallel with Section \ref{section06}, we complete the algebra $\tYht$ as follows. Let $I_p$ with $p\geqslant 0$ be the two-sided ideal in  $\tYht$ generated by $h^{p+1}$ and all $\tau_{ij}^{(-r)}$ with $r>p$. Then the completed algebra is defined by
$$
\tYhtc =\lim_{\longleftarrow} \tYht/I_p.
$$
 By arguing as in the previous section, one obtains from 
\eqref{trigbethertt} the identity
$$
   \Tc_1(x)\ts \Tc_2 (y)
=  \wht{\Rc} ( x/y )^{-1}  \Tc_2 (y)\ts \Tc_1(x)\ts\wht{\Rc} ( x/y  ),
$$
where the coefficients of its matrix entries   are well-defined elements of    $\tYhtc$; recall \eqref{bethecomp}.

For any   $n=1,\ldots ,N$ introduce the series 
$$
\Tc_{(n)}  (z)=
\tr_{1,\ldots ,n}\ts A^h_{(n)} \ts\Tc_1(z )\Tc_2(ze^{-h})\ldots \Tc_n(ze^{-(n-1)h})\in \tYhtc[[z^{\pm 1}]].
$$
Proposition \ref{bethemc}  is generalized to the trigonometric setting as follows.

\begin{pro}\label{trigbethemc}
The coefficients of  all $\Tc_{(n)}  (z)$
  mutually commute.
\end{pro}

\begin{prf}
(1) As with   Proposition \ref{bethemc}, we start by establishing some preliminary identities which are needed in the main part of the proof.
The fusion procedure \eqref{trigfusion} implies    that by applying the substitutions 
$x_i=ze^{-(i-1)h}$, $i=1,\ldots ,n$
to the expression
$$
\Tc_{[n]}(x)\coloneqq \left(\prod_{i=1,\dots,n-1}^{\longrightarrow} 
\prod_{j=i+1,\ldots,n}^{\longrightarrow}  \wht{\Rc}_{ij}(  x_i /x_j )\right) \Tc_1(x_1)\ldots \Tc_n (x_n) 
$$
we obtain, up to a nonzero scalar multiple,
$$
\Tc_{[n]}(z,ze^{-h},\ldots, ze^{-(n-1)h})=A^h_{(n)} \ts\Tc_1(z )\Tc_2(ze^{-h})\ldots \Tc_n(ze^{-(n-1)h})
.
$$

The   notation conventions for     $R$-matrix products  from  Section \ref{section01}   naturally generalize  to the trigonometric case, e.g., for $u=(u_1,\ldots ,u_n)$ and $v=(v_1,\ldots ,v_m)$ we have (cf. \eqref{oppositeof2f})
$$
\wht{\Rc}_{nm}^{12}(u/v)= \prod_{i=1,\dots,n}^{\longrightarrow} 
\prod_{j=n+1,\ldots,n+m}^{\longleftarrow}   \wht{\Rc}_{ij} ( u_i/v_{j-n}  ).
$$
Other consequences of \eqref{trigfusion} to be needed   for the main part of the proof  are  
\begin{gather}
A^{h\ts 1}_{(n_1 )} \ts A^{h\ts 2}_{(n_2)} \ts \wht{\Rc}_{n_1 n_2}^{12}(x_1 /x_2  )
= \wht{\Rc}_{n_1 n_2}^{\bar{\scriptstyle 1}\bar{\scriptstyle 2}}( x_1 /x_2 ) \ts
A^{h\ts 1}_{(n_1 )} \ts A^{h\ts 2}_{(n_2)} ,\non\\
 \Tc_{[n_1]}^{13}(x_1)\ts \Tc_{[n_2]}^{23} (x_2)
=
\wht{\Rc}_{n_1 n_2}^{\bar{\scriptstyle 1}\bar{\scriptstyle 2}}( x_1/x_2 )^{-1}
\Tc_{[n_2]}^{23} (x_2)\ts \Tc_{[n_1]}^{13}(x_1)\ts \wht{\Rc}_{n_1 n_2}^{ 12}(  x_1/x_2   ),\label{t6z7u}\\
A^h_{(n_i) } \ts  \Tc_1(z_i) \ts  \Tc_2(z_i e^{-h})\ldots \Tc_{n_i}(z_i e^{-(n_i-1)h}) = \Tc_{n_i}(z_i e^{-(n_i-1)h})\ldots \Tc_2(z_i e^{-h})\ts \Tc_1(z_i)\ts A^h_{(n_i) },\label{t6z7v}
\end{gather}
 where $A^{h\ts 1}_{(n_1 )}$ (resp. $A^{h\ts 2}_{(n_2)}$) denotes the action of the $h$-anti-symmetrizer on the tensor factors $1,\ldots ,n_1$ (resp. $n_1+1,\ldots, n_1+n_2$), $z_1$ and $z_2$ are single variables and 
$$x_i=(z_i,z_i e^{-h},\ldots ,z_i e^{-(n_i-1)h})\quad\text{for }i=1,2.$$
They can be easily  verified by using the Yang--Baxter equation \eqref{trigybe} and   defining relations \eqref{trigbethertt}.
Note that their rational counterparts  are given by \eqref{betheeyb2}, \eqref{thrltn} and \eqref{betheeyb3}.\vspace{2pt}

\noindent (2) The  proposition   can be now proved by repeating the arguments from the second part of the proof of Proposition \ref{bethemc} and employing the identities listed above.  
\end{prf}

Finally, consider the series
$$
\mathbb{T}(z) \coloneqq \mathcal{T}_{(N)}(z)=
\tr_{1,\ldots ,N}\ts A^h_{(N)}\ts\Tc_1(z )\Tc_2(ze^{-h})\ldots \Tc_N(ze^{-(N-1)h}) .
$$
By using   property   \eqref{qantisimprop} of the  anti-symmetrizer
and  arguing as in the proof of Theorem \ref{centralnithm}, one obtains its trigonometric counterpart:
\begin{thm}\label{trigcentralnithm}
The coefficients of    $\mathbb{T}(z)$
belong to the center of $\tYhtc$.
\end{thm}

 \begin{rem}
The results from Sections \ref{sec0201}--\ref{sec07}, along with Remark \ref{remarkrmatricca},   
 indicate a possible interpretation of Proposition \ref{trigbethemc} and Theorem \ref{trigcentralnithm}
from the viewpoint of theory of $\phi$-coordinated modules for quantum vertex algebras; cf. \cite[Sect. 3]{JKLT} and \cite{Liphi}. 
  \end{rem}

\section*{Acknowledgement} This work has been supported in part by Croatian Science Foundation under the project UIP-2019-04-8488.

\linespread{1.0}

\end{document}